# ALGEBRAIC METHODS TOWARD HIGHER-ORDER PROBABILITY INEQUALITIES, II

By Donald St. P. Richards

*Pennsylvania State University*

Let $(L, \preccurlyeq)$ be a finite distributive lattice, and suppose that the functions $f_1, f_2 : L \to \mathbb{R}$ are monotone increasing with respect to the partial order $\preccurlyeq$. Given $\mu$ a probability measure on $L$, denote by $\mathbb{E}(f_i)$ the average of $f_i$ over $L$ with respect to $\mu$, $i = 1, 2$. Then the FKG inequality provides a condition on the measure $\mu$ under which the covariance, $\mathrm{Cov}(f_1, f_2) := \mathbb{E}(f_1 f_2) - \mathbb{E}(f_1)\mathbb{E}(f_2)$, is nonnegative. In this paper we derive a "third-order" generalization of the FKG inequality: Let $f_1$, $f_2$ and $f_3$ be nonnegative, monotone increasing functions on $L$; and let $\mu$ be a probability measure satisfying the same hypotheses as in the classical FKG inequality; then

$$2\mathbb{E}(f_1 f_2 f_3)$$
$$- [\mathbb{E}(f_1 f_2)\mathbb{E}(f_3) + \mathbb{E}(f_1 f_3)\mathbb{E}(f_2) + \mathbb{E}(f_1)\mathbb{E}(f_2 f_3)]$$
$$+ \mathbb{E}(f_1)\mathbb{E}(f_2)\mathbb{E}(f_3)$$

is nonnegative. This result reduces to the FKG inequality for the case in which $f_3 \equiv 1$.

We also establish fourth- and fifth-order generalizations of the FKG inequality and formulate a conjecture for a general $m$th-order generalization. For functions and measures on $\mathbb{R}^n$ we establish these inequalities by extending the method of diffusion processes. We provide several applications of the third-order inequality, generalizing earlier applications of the FKG inequality. Finally, we remark on some connections between the theory of total positivity and the existence of inequalities of FKG-type within the context of Riemannian manifolds.

**1. Introduction.** In the realm of probability inequalities, the FKG inequality, due to Fortuin, Kasteleyn and Ginibre (1971), now occupies a po-sition of fundamental importance because of its simplicity and widespread applicability. Before we state in detail this remarkable inequality, it is worthwhile to note some of its many applications.

In statistics, the FKG inequality has appeared in the study of monotonicity properties of power functions of likelihood ratio test statistics in multivariate analysis, association and dependence properties of random variables, and observational studies; see Sarkar (1969), Perlman and Olkin (1980), Eaton (1987), Karlin and Rinott (1980) and Rosenbaum (1995).

In probability theory and mathematical physics, the FKG inequality has appeared in the areas of diffusion equations, interacting particle systems, Ising models, reliability theory and percolation; see Grimmett (1999), Herbst and Pitt (1991), Lebowitz (1972), Liggett (1985), Preston (1976), Simon (1974, 1979) and Glimm and Jaffe (1987).

In work on total positivity and its connections with the theory of finite reflection groups and analysis on Lie groups, an analog of the FKG inequality was derived by Gross and Richards (1995).

In combinatorial theory, the FKG inequality has appeared in work on the monotonicity of partial orders, Sperner theory, graph theory and Ramsey theory; see Graham (1982, 1983), Engel (1997), Cameron (1987) and Seymour and Welsh (1975).

We now state the inequality. Let $L$ be a finite distributive lattice with partial ordering $\preccurlyeq$, least upper bound $\vee$ and greatest lower bound $\wedge$. A function $f: L \to \mathbb{R}$ is called (*monotone*) *increasing* if $f(x) \leq f(y)$ whenever $x \preccurlyeq y$. A probability measure $\mu$ on $L$ is said to be *multivariate totally positive of order* 2 (MTP$_2$) if

$$\mu(x \vee y)\mu(x \wedge y) \geq \mu(x)\mu(y) \tag{1.1}$$

for all $x, y \in L$. In some parts of the literature, an MTP$_2$ probability measure $\mu$ is called *an FKG measure* or *log-supermodular*; we prefer the MTP$_2$ terminology, reflecting the relationship with the classical theory of total positivity wherein (1.1) is an abstract formulation of the concept of total positivity of order 2.

For any probability measure $\mu$ on $L$ and any function $f: L \to \mathbb{R}$, denote by

$$\mathbb{E}(f) := \sum_{a \in L} \mu(a) f(a)$$

the *mean* or *average* of $f$ with respect to $\mu$.

Suppose $f_1$ and $f_2$ are both increasing (or both decreasing) real-valued functions on $L$, and let $\mu$ be an MTP$_2$ probability measure on $L$. Then the FKG inequality provides that

$$\mathrm{Cov}(f_1, f_2) := \mathbb{E}(f_1 f_2) - \mathbb{E}(f_1)\mathbb{E}(f_2) \geq 0. \tag{1.2}$$

In mathematical statistics, it is usual to state the FKG inequality for the space $\mathbb{R}^n$. In that setting, for vectors $x = (x_1, \ldots, x_n)$ and $y = (y_1, \ldots, y_n)$ in $\mathbb{R}^n$, the corresponding partial order is given by $x \preccurlyeq y$ if $x_j \leq y_j$ for all $j = 1, \ldots, n$; the lattice operations $\vee$ and $\wedge$ are

$$x \vee y = (\max(x_1, y_1), \ldots, \max(x_n, y_n))$$

and

$$x \wedge y = (\min(x_1, y_1), \ldots, \min(x_n, y_n));$$

a function $f : \mathbb{R}^n \to \mathbb{R}$ is *increasing* if $f(x) \leq f(y)$ whenever $x \preccurlyeq y$; and a probability density function $K : \mathbb{R}^n \to \mathbb{R}$, or its underlying random vector $(X_1, \ldots, X_n)$, is said to be *multivariate totally positive of order* 2 ($\text{MTP}_2$) if

$$K(x \vee y) K(x \wedge y) \geq K(x) K(y) \tag{1.3}$$

for all $x, y \in \mathbb{R}^n$. If the functions $f_1, f_2 : \mathbb{R}^n \to \mathbb{R}$ are both increasing or both decreasing and the expectations $E(f_1)$, $E(f_2)$ and $E(f_1 f_2)$, with respect to $K$, are finite, then the FKG inequality on $\mathbb{R}^n$ provides that

$$\begin{aligned}
&\int_{\mathbb{R}^n} f_1(x) f_2(x) K(x) \, dx \\
&\quad - \left( \int_{\mathbb{R}^n} f_1(x) K(x) \, dx \right) \cdot \left( \int_{\mathbb{R}^n} f_2(x) K(x) \, dx \right) \geq 0.
\end{aligned} \tag{1.4}$$

It is well known that the continuous case of the FKG inequality, (1.4), can be deduced by an approximation argument from the discrete case, (1.2); see Karlin and Rinott (1980). We also remark that numerous generalizations of (1.2) have appeared in the literature; see Ahlswede and Daykin (1978), Batty and Bollmann (1980), Holley (1974), Preston (1974), Kemperman (1977), Edwards (1978) and Rinott and Saks (1993). These results typically provide conditions leading to inequalities of the form

$$\int \prod_{i=1}^k f_i(x) g_i(x) \cdot K(x) \, dx \\
- \left( \int \prod_{i=1}^k f_i(x) \cdot K(x) \, dx \right) \left( \int \prod_{i=1}^k g_i(x) \cdot K(x) \, dx \right) \geq 0$$

for classes of nonnegative functions $f_1, \ldots, f_k$ and $g_1, \ldots, g_k$.

In this paper we derive generalizations of (1.2) or (1.4) involving alternating sums with more than two terms. Recall that the covariance between two random variables is an example of a cumulant (or Ursell function) of those random variables, so the FKG inequality (1.2) provides an inequality for the simplest cumulant of the random variables $f_1$ and $f_2$. Hence, in a search for "higher-order" generalizations of (1.2) involving several functions $f_1, f_2, \ldots$, it is natural that we study the higher cumulants.

There are well-known probability distributions for which the FKG inequality holds but for which the higher cumulants are nonpositive. Indeed, in the case of any Gaussian distribution, the higher cumulants are identically zero. It is then apparent that, without additional restrictions on the density function $K$, the cumulants themselves cannot provide direct generalizations of (1.2) [in this regard, we refer to Percus (1975) and Sylvester (1975) for other types of correlation inequalities for Ursell functions].

Nevertheless, the algebraic structure of the higher-order cumulants provides crucial motivation for generalizing the FKG inequality. To explain, recall that the alternating sum

$$
\begin{aligned}
\kappa_3(f_1, f_2, f_3) :=\ & \mathbb{E}(f_1 f_2 f_3) \\
& - [\mathbb{E}(f_1 f_2)\mathbb{E}(f_3) + \mathbb{E}(f_1 f_3)\mathbb{E}(f_2) + \mathbb{E}(f_1)\mathbb{E}(f_2 f_3)] \\
& + 2\mathbb{E}(f_1)\mathbb{E}(f_2)\mathbb{E}(f_3)
\end{aligned}
\tag{1.5}
$$

is a cumulant of the random variables $f_1$, $f_2$ and $f_3$ [cf. Speed (1983)]. We define the *third-order conjugate cumulant*

$$
\begin{aligned}
\kappa_3'(f_1, f_2, f_3) :=\ & 2\mathbb{E}(f_1 f_2 f_3) \\
& - [\mathbb{E}(f_1 f_2)\mathbb{E}(f_3) + \mathbb{E}(f_1 f_3)\mathbb{E}(f_2) + \mathbb{E}(f_1)\mathbb{E}(f_2 f_3)] \\
& + \mathbb{E}(f_1)\mathbb{E}(f_2)\mathbb{E}(f_3),
\end{aligned}
\tag{1.6}
$$

which is derived from $\kappa_3$ by "reversing" the order of the absolute value of the coefficients appearing in (1.5). The general $m$th-order conjugate cumulant will be defined in a similar manner. Our first main result, proven in Section 2, is the following.

THEOREM 1.1. *Let $L$ be a finite distributive lattice, let $\mu$ be an $MTP_2$ probability measure on $L$, and let $f_1$, $f_2$ and $f_3$ be nonnegative increasing functions on $L$. Then*

$$\kappa_3'(f_1, f_2, f_3) \geq 0. \tag{1.7}$$

As a consequence of Theorem 1.1, we obtain a lower bound on the cumulant $\kappa_3(f_1, f_2, f_3)$.

COROLLARY 1.2. *Under the same hypotheses as Theorem 1.1,*

$$\kappa_3(f_1, f_2, f_3) \geq -[\mathbb{E}(f_1 f_2 f_3) - \mathbb{E}(f_1)\mathbb{E}(f_2)\mathbb{E}(f_3)].$$

That the inequality (1.7) generalizes the FKG inequality can be seen in three ways. First, (1.7) reduces to (1.2) if $f_3 \equiv 1$. Second, (1.7) does not generally reduce to (1.2) if, say, $f_3$ is the indicator function of a proper subset of $L$. Third, it is straightforward to verify that (1.7) can be rewritten as

$$\mathrm{Cov}(f_1 f_2, f_3) - \mathbb{E}(f_1)\mathrm{Cov}(f_2, f_3) + \mathrm{Cov}(f_1 f_3, f_2) \geq 0,$$

expressing the property that an alternating linear combination of nonnegative covariances is nonnegative. It is also a pleasant surprise that these results hold under the same hypotheses on $K$ required for the classical FKG inequality, so that the only additional assumptions required are the existence of the various expectations.

With Theorem 1.1 in place, it is natural to ask for similar generalizations of the FKG inequality involving four or more functions. To that end, we shall prove the following results.

THEOREM 1.3. *Let $\mu$ be an $MTP_2$ probability measure on $L$ and let $f_i$, $i = 1, \ldots, 5$, be nonnegative increasing functions on $L$. Then the fourth- and fifth-order conjugate cumulants,*

$$\begin{aligned}\kappa_4' := {} & 6\mathbb{E}(f_1 f_2 f_3 f_4) \\ & - 2[\mathbb{E}(f_1 f_2 f_3)\mathbb{E}(f_4) + \cdots] - [\mathbb{E}(f_1 f_2)\mathbb{E}(f_3 f_4) + \cdots] \\ & + [\mathbb{E}(f_1 f_2)\mathbb{E}(f_3)\mathbb{E}(f_4) + \cdots] - \mathbb{E}(f_1)\mathbb{E}(f_2)\mathbb{E}(f_3)\mathbb{E}(f_4)\end{aligned} \quad (1.8)$$

*and*

$$\begin{aligned}\kappa_5' := {} & 24\mathbb{E}(f_1 f_2 f_3 f_4 f_5) \\ & - 6[\mathbb{E}(f_1 f_2 f_3 f_4)\mathbb{E}(f_5) + \cdots] - 2[\mathbb{E}(f_1 f_2 f_3)\mathbb{E}(f_4 f_5) + \cdots] \\ & + 2[\mathbb{E}(f_1 f_2 f_3)\mathbb{E}(f_4)\mathbb{E}(f_5) + \cdots] + [\mathbb{E}(f_1 f_2)\mathbb{E}(f_3 f_4)\mathbb{E}(f_5) + \cdots] \\ & - [\mathbb{E}(f_1 f_2)\mathbb{E}(f_3)\mathbb{E}(f_4)\mathbb{E}(f_5) + \cdots] + \mathbb{E}(f_1)\mathbb{E}(f_2)\mathbb{E}(f_3)\mathbb{E}(f_4),\end{aligned} \quad (1.9)$$

*are nonnegative.*

In the above results, the notation, for example, "$\mathbb{E}(f_1 f_2)\mathbb{E}(f_3)\mathbb{E}(f_4) + \cdots$" is shorthand notation for the sum over all distinct terms consisting of products of expectations of the set of functions $\{f_1, \ldots, f_4\}$, divided into subsets of sizes 2, 1, and 1. Explicitly,

$$\begin{aligned}& \mathbb{E}(f_1 f_2)\mathbb{E}(f_3)\mathbb{E}(f_4) + \cdots \\ & \quad \equiv \mathbb{E}(f_1 f_2)\mathbb{E}(f_3)\mathbb{E}(f_4) + \mathbb{E}(f_1 f_3)\mathbb{E}(f_2)\mathbb{E}(f_4) + \mathbb{E}(f_1 f_4)\mathbb{E}(f_2)\mathbb{E}(f_3) \\ & \qquad + \mathbb{E}(f_1)\mathbb{E}(f_2 f_3)\mathbb{E}(f_4) + \mathbb{E}(f_1)\mathbb{E}(f_2 f_4)\mathbb{E}(f_3) + \mathbb{E}(f_1)\mathbb{E}(f_2)\mathbb{E}(f_3 f_4).\end{aligned}$$

In Section 2 we shall prove Theorem 1.1. Once this is complete, the proof of Theorem 1.3 is seen to be a consequence of elementary, but lengthy, algebraic calculations which are easily performed using computer algebra software, for example, MAPLE. Further, we shall conjecture a generalization of Theorems 1.1 and 1.3 to the case of an arbitrary number of functions.

In Section 3 we apply the method of diffusion processes to provide another approach to the generalized FKG inequalities. We extend a method due to Herbst and Pitt (1991) who derived the FKG inequality in their work on diffusion semigroups which are stochastically monotone and which preserve the class of positively correlated measures on $\mathbb{R}^n$. While the methods utilized in Section 2 are patently algebraic in flavor, the methods used in Section 3 may initially appear to be analytic in nature; however, a deeper reading will also reveal the ideas of this latter section to be as algebraic in characteristic as are those of Section 2.

In Section 4 we approach the problem of proving Theorem 1.1 through the method of duplicate variables; this approach seems significantly more complex than in the case of the classical FKG inequality, and we provide explicit details for the case in which $n = 1$. In Section 5 we provide a number of applications, extending some well-known applications of the classical FKG inequality. Finally, in Section 6 we remark on connections between the theory of total positivity and the possible existence of inequalities of FKG-type within the context of weakly symmetric Riemannian manifolds.

**2. Finite distributive lattices.** Let $L$ be a finite distributive lattice with partial order $\preccurlyeq$, least upper bound $\vee$ and greatest lower bound $\wedge$. Recall that every finite distributive lattice is order-isomorphic to the lattice of subsets of a finite set. Then there exists a finite set $A$ such that $L = 2^A$, the collection of all subsets of $A$. Endowed with set inclusion as the partial order, and with union and intersection as least upper bound and greatest lower bound, respectively, $2^A$ is a finite distributive lattice.

A probability measure $\mu$ on $2^A$ is said to be *multivariate totally positive of order* 2 (MTP$_2$) if, for all $a, b \subseteq A$,

$$\mu(a \cup b)\mu(a \cap b) \geq \mu(a)\mu(b). \tag{2.1}$$

A function $f : 2^A \to \mathbb{R}$ is called (*monotone*) *increasing* if $f(a) \geq f(b)$ whenever $a, b \subseteq A$ satisfy $a \supseteq b$. We denote the *expected value* of $f$ with respect to the measure $\mu$ by

$$\mathbb{E}(f) := \sum_{a \subseteq A} \mu(a) f(a).$$

We now establish (1.7).

PROOF OF THEOREM 1.1. Our argument follows that of den Hollander and Keane (1986). Without loss of generality, we assume that $L = 2^A$ for some finite set $A$. If the cardinality of $2^A$ is 1, that is, $A = \varnothing$, then (1.7) holds trivially; therefore we may suppose that $A$ is nonempty.

Let us first assume that $\mu(a) > 0$ for all $a \subseteq A$. Choose and fix $B$, an arbitrarily chosen subset of $A$. For any $a \subseteq B$, define

$$\mu_B(a) := \sum_{b \subseteq A \setminus B} \mu(a \cup b) \tag{2.2}$$

and, for any function $f : 2^A \to \mathbb{R}$, define

$$f_B(a) := \frac{1}{\mu_B(a)} \sum_{b \subseteq A \setminus B} \mu(a \cup b) f(a \cup b). \tag{2.3}$$

As observed by den Hollander and Keane (1986), it is clear that $\mu_B$ is the marginal probability measure on the lattice $2^B$; also, $f_B(a)$ is the conditional expectation of $f$, given $a \subseteq B$. It can be shown [cf. den Hollander and Keane (1986), page 171] that $\mu_B$ is MTP$_2$, and that $f_B$ is increasing if $f$ is increasing [see Eaton (1987) or Karlin and Rinott (1980) for similar results in the case of $\mathbb{R}^n$].

For any function $g : 2^B \to \mathbb{R}$, we define

$$\mathbb{E}_B(g) := \sum_{a \subseteq B} \mu_B(a) g(a).$$

In the sequel we shall need the double expectation theorem: For any $B \subseteq A$,

$$\mathbb{E}(f) = \sum_{a \subseteq A} \mu(a) f(a) = \sum_{a \subseteq B} \mu_B(a) f_B(a) = \mathbb{E}_B(f_B). \tag{2.4}$$

In words, the expected value of $f$ equals the expected value of its conditional expectations.

To establish (1.7) it suffices, as observed by den Hollander and Keane (1986) in the case of the FKG inequality, to assume that $B = A \setminus \{z\}$, where $z \in A$ is chosen arbitrarily; this amounts to a proof by induction on the length of maximal chains in the lattice $2^A$. Using the shorthand notation $f_{iB}$ to denote $(f_i)_B$, $i = 1, 2, 3$, we claim that

$$\begin{aligned}(2.5) \quad & 2\mathbb{E}_B((f_1 f_2 f_3)_B) \\ & - [\mathbb{E}_B((f_1 f_2)_B f_{3B}) + \mathbb{E}_B((f_1 f_3)_B f_{2B}) + \mathbb{E}_B(f_{1B}(f_2 f_3)_B)] \\ & + \mathbb{E}_B(f_{1B} f_{2B} f_{3B}) \geq 0.\end{aligned}$$

By (2.2) and (2.3) we have

$$\mu_B(a) = \mu(a) + \mu(a \cup \{z\}) \tag{2.6}$$

and

$$f_B(a) = \frac{1}{\mu_B(a)} (\mu(a) f(a) + \mu(a \cup \{z\}) f(a \cup \{z\})) \tag{2.7}$$

for $a \subseteq B$. By (2.7) we have

$$\mu_B(a)^3 (f_1 f_2 f_3)_B(a)$$
$$(2.8) \qquad = \mu_B(a)^2 [\mu(a) f_1(a) f_2(a) f_3(a)$$
$$\qquad\qquad + \mu(a \cup \{z\}) f_1(a \cup \{z\}) f_2(a \cup \{z\}) f_3(a \cup \{z\})].$$

Further, for $\{i, j, k\} = \{1, 2, 3\}$,

$$\mu_B(a)^3 (f_i f_j)_B(a) f_{kB}(a)$$
$$(2.9) \qquad = \mu_B(a)[\mu(a) f_i(a) f_j(a) + \mu(a \cup \{z\}) f_i(a \cup \{z\}) f_j(a \cup \{z\})]$$
$$\qquad\qquad \times [\mu(a) f_k(a) + \mu(a \cup \{z\}) f_k(a \cup \{z\})]$$

and

$$(2.10) \qquad \mu_B(a)^3 f_{1B}(a) f_{2B}(a) f_{3B}(a)$$
$$= \prod_{i=1}^{3} [\mu(a) f_i(a) + \mu(a \cup \{z\}) f_i(a \cup \{z\})].$$

Collecting together (2.8)–(2.10) and performing elementary algebraic simplifications with the aid of (2.6) and (2.7), we obtain

$$\mu_B(a)^3 \{2 (f_1 f_2 f_3)_B(a)$$
$$\qquad - [(f_1 f_2)_B(a) f_{3B}(a) + (f_1 f_3)_B(a) f_{2B}(a) + f_{1B}(a)(f_2 f_3)_B(a)]$$
$$(2.11) \qquad\qquad\qquad\qquad\qquad + f_{1B}(a) f_{2B}(a) f_{3B}(a)\}$$
$$= \mu(a) \mu(a \cup \{z\})$$
$$\qquad \times [(f_1(a \cup \{z\}) - f_1(a)) \Phi_{1B}(a) + f_1(a \cup \{z\}) \Phi_{2B}(a)],$$

where

$$\Phi_{1B}(a) = \mu(a) (f_2(a \cup \{z\}) f_3(a \cup \{z\}) - f_2(a) f_3(a))$$
$$(2.12) \qquad\qquad + \mu(a \cup \{z\}) f_3(a) (f_2(a \cup \{z\}) - f_2(a))$$
$$\qquad\qquad + \mu(a \cup \{z\}) f_2(a) (f_3(a \cup \{z\}) - f_3(a))$$

and

$$\Phi_{2B}(a) = (\mu(a \cup \{z\}) + \mu(a))$$
$$(2.13) \qquad\qquad \times (f_2(a \cup \{z\}) - f_2(a))(f_3(a \cup \{z\}) - f_3(a)).$$

Since each $f_i$ is nonnegative and increasing, it follows that (2.12) and (2.13) are sums of products of nonnegative terms; hence (2.11) is nonnegative.

Next we divide both sides of (2.11) by $\mu_B(a)^2$ and sum over all $a \subseteq B$. We have

$$\sum_{a \subseteq B} \mu_B(a)(f_1 f_2 f_3)_B(a) = \mathbb{E}_B((f_1 f_2 f_3)_B) = \mathbb{E}(f_1 f_2 f_3),$$

the latter equality following from the double expectation theorem (2.4). For $\{i, j, k\} = \{1, 2, 3\}$, we have

$$\sum_{a \subseteq B} \mu_B(a)(f_i f_j)_B(a) f_{kB}(a) = \mathbb{E}_B((f_i f_j)_B f_{kB});$$

and also

$$\sum_{a \subseteq B} \mu_B(a) f_{1B}(a) f_{2B}(a) f_{3B}(a) = \mathbb{E}_B(f_{1B} f_{2B} f_{3B}).$$

Collecting these identities together and applying the nonnegativity of (2.11), we obtain (2.5).

Since $B$ was chosen arbitrarily, we may set $B = \varnothing$ in (2.5). By (2.2) and (2.3) we have $\mathbb{E}_\varnothing(f_{i\varnothing}) \equiv \mathbb{E}(f_i)$, $\mathbb{E}_\varnothing((f_i f_j)_\varnothing f_{k\varnothing}) \equiv \mathbb{E}(f_i f_j)\mathbb{E}(f_k)$, where $\{i, j, k\} = \{1, 2, 3\}$, and also $\mathbb{E}_\varnothing(f_{1\varnothing} f_{2\varnothing} f_{3\varnothing}) \equiv \mathbb{E}(f_1)\mathbb{E}(f_2)\mathbb{E}(f_3)$. Then (2.5) reduces to (1.7).

Finally, the case in which $\mu$ is not everywhere positive is resolved as in the case of the FKG inequality [see den Hollander and Keane (1986), page 172]. That is, we can carry through the above arguments once $\mu_B$ and the $f_{iB}$ are defined on the set $\{a \subseteq B : \mu_B(a) \neq 0\}$; this is done in any way which ensures that all $f_{iB}$ are increasing on $2^B$, and the actual choice is immaterial since all we need is that each $f_{iB}$ is increasing on the support of $\mu_B$. □

REMARK 2.1. (i) The assumption that the functions $f_1$, $f_2$ and $f_3$ are nonnegative is essential; in particular, we cannot avoid this assumption by adding a constant to each function, for $\kappa_3'(f_1, f_2, f_3)$ is not invariant under such shifts. A referee has noted that without the positivity assumption, $E(X^3)$ can be made arbitrarily negatively large relative to other moments, and this would violate the inequality for $f_i(x) = x$, $i = 1, 2, 3$.

(ii) Consider the range of values of positive coefficients $c_1$, $c_2$ and $c_3$ such that

$$c_1 \mathbb{E}(f_1 f_2 f_3) - c_2[\mathbb{E}(f_1 f_2)\mathbb{E}(f_3) + \mathbb{E}(f_1 f_3)\mathbb{E}(f_2) + \mathbb{E}(f_1)\mathbb{E}(f_2 f_3)]$$
$$+ c_3 \mathbb{E}(f_1)\mathbb{E}(f_2)\mathbb{E}(f_3) \geq 0$$

for all nonnegative increasing functions $f_1$, $f_2$ and $f_3$. For simplicity, let us work with functions on $\mathbb{R}^2$. In order to make comparison with Theorem 1.1, we obviously need to impose the restriction

$$c_1 - 3c_2 + c_3 = 0.$$

Without loss of generality, we may assume $c_2 = 1$. By setting $f_3 \equiv 1$ we obtain

$$(c_1 - 1)\mathbb{E}(f_1 f_2) - (2 - c_3)\mathbb{E}(f_1)\mathbb{E}(f_2) \geq 0.$$

In order to maintain positive coefficients, we must impose the additional restriction $c_1 \geq 1$.

For $j = 1, 2, 3$, let $a_j, b_j \in \mathbb{R}$ and let $f_j$ be the indicator function of the set $\{(u, v) \in \mathbb{R}^2 : u \geq a_j, v \geq b_j\}$. Denote the underlying random vector by $(X, Y)$; by choosing suitable distributions for $(X, Y)$ it can be shown that the condition $c_1 \geq 1$ is not sufficient to ensure that $\kappa'_3(f_1, f_2, f_3) \geq 0$. For instance, suppose $a_1 \leq a_2 \leq a_3$ and $b_1 \leq b_2 \leq b_3$. Denoting $P(X \geq a_j, Y \geq b_j)$ by $\pi_j$, we obtain

$$\kappa'_3(f_1, f_2, f_3) = c_1 \pi_3 - (2\pi_2 \pi_3 + \pi_1 \pi_3) + c_3 \pi_1 \pi_2 \pi_3$$
$$\equiv \pi_3[(1 - \pi_2)(2 + (c_1 - 3)\pi_1) + (c_1 - 2)(1 - \pi_1)].$$

For $c_1 \geq 1$ we then have

$$\kappa'_3(f_1, f_2, f_3) \geq \pi_3(1 - \pi_1)(1 - 2\pi_2).$$

For suitable $(X, Y)$ it is possible to attain this lower bound; moreover, the bound can be negative. Therefore the restriction $c_1 \geq 1$ is not sufficient to ensure nonnegativity of $\kappa'_3(f_1, f_2, f_3)$ for all nonnegative increasing $f_j$. By choosing the $f_j$ from among the class of indicator functions of the type above, we can deduce that $c_1 \geq 2$ is sufficient. We will return to this theme in Section 6.

Returning to the general context, a close inspection of the proof of Theorem 1.1 shows that we have also obtained a collection of lower bounds for $\mathbb{E}(f_1 f_2 f_3)$.

COROLLARY 2.2. *Let $\mu$ be an $MTP_2$ measure on $2^A$ and let $f_1$, $f_2$ and $f_3$ be nonnegative increasing functions on $2^A$. For any $B \subseteq A$, there holds the lower bound*

$$\begin{aligned}(2.14) \quad 2\mathbb{E}(f_1 f_2 f_3) &\geq \mathbb{E}_B((f_1 f_2)_B f_{3B}) + \mathbb{E}_B((f_1 f_3)_B f_{2B}) \\ &\quad + \mathbb{E}_B(f_{1B}(f_2 f_3)_B) - \mathbb{E}_B(f_{1B} f_{2B} f_{3B}).\end{aligned}$$

REMARK 2.3. In the case of the FKG inequality, that is, for the case in which $f_3 \equiv 1$, the left-hand side of (2.5) reduces to

$$(2.15) \quad \mathbb{E}_B((f_1 f_2)_B) - \mathbb{E}_B(f_{1B} f_{2B}) \equiv \mathbb{E}(f_1 f_2) - \mathbb{E}_B(f_{1B} f_{2B}).$$

In addition to establishing the FKG inequality, den Hollander and Keane [(1986), Theorem 4(b)] establish the sharper inequality that (2.15) is a decreasing function of $B \in 2^A$. This raises the issue of whether the left-hand side of (2.5) satisfies similar monotonicity properties. Even for the simplest lattices, such monotonicity properties appear difficult to discern. In particular, as the following counterexample shows, the left-hand side of (2.5) is not generally monotonically decreasing in $B$.

Suppose that $A = \{w, z\}$, a set with two distinct elements. The corresponding lattice is $2^A = \{\varnothing, \{w\}, \{z\}, \{w, z\}\}$. Define a probability measure $\mu$ on $2^A$ by: $\mu(\phi) = 1/2$, $\mu(\{w\}) = \mu(\{z\}) = 1/8$ and $\mu(\{w,z\}) = 1/4$; then it is straightforward to verify that $\mu$ is MTP$_2$. Define three nonnegative increasing functions $f_1$, $f_2$ and $f_3$ on $2^A$ by the substitutions

$$f_i(a) = \begin{cases} \alpha_i, & \text{if } a = \varnothing, \\ \alpha_i + \beta_i, & \text{if } a = \{w\}, \\ \alpha_i + \gamma_i, & \text{if } a = \{z\}, \\ \alpha_i + \beta_i + \gamma_i + \delta_i, & \text{if } a = \{w, z\}, \end{cases}$$

where $\alpha_i, \beta_i, \gamma_i, \delta_i \geq 0$, $i = 1, 2, 3$. If we denote by $g(B)$ the left-hand side of (2.5), then for the case in which $(\alpha_1, \alpha_2, \alpha_3) = (\beta_1, \beta_2, \beta_3) = (1, 2, 3)$, $(\gamma_1, \gamma_2, \gamma_3) = (4, 5, 6)$ and $(\delta_1, \delta_2, \delta_3) = (0.1, 0.2, 0.3)$, a straightforward calculation (carried out using MAPLE) reveals that $g(\varnothing) - g(\{w\}) < 0$. On the other hand, for the case in which $\alpha_i = \beta_i = \gamma_i = \delta_i = i$, $i = 1, 2, 3$, we have $g(\varnothing) - g(\{w\}) > 0$.

REMARK 2.4. Before we turn to the proof of Theorem 1.3, let us review the algebraic calculations appearing in the proof of Theorem 1.1. In (2.6) and (2.7), introduce the notation $u_i = f_i(a \cup \{z\})$, $v_i = f_i(a)$, $i = 1, 2, 3$; and let $\omega_1 = \mu(a \cup \{z\})$ and $\omega_2 = \mu(a)$ be nonnegative weights. With $u = (u_1, u_2, u_3)$ and $v = (v_1, v_2, v_3)$ as variables, define the polynomials

$$p_{(3)}(u; v) := (\omega_1 + \omega_2)^2 (\omega_1 u_1 u_2 u_3 + \omega_2 v_1 v_2 v_3),$$

$$p_{(2,1)}(u; v) := (\omega_1 + \omega_2)[(\omega_1 u_1 u_2 + \omega_2 v_1 v_2)(\omega_1 u_3 + \omega_2 v_3)$$
$$+ (\omega_1 u_1 u_3 + \omega_2 v_1 v_3)(\omega_1 u_2 + \omega_2 v_2)$$
$$+ (\omega_1 u_2 u_3 + \omega_2 v_2 v_3)(\omega_1 u_1 + \omega_2 v_1)],$$

$$p_{(1,1,1)}(u; v) := (\omega_1 u_1 + \omega_2 v_1)(\omega_1 u_2 + \omega_2 v_2)(\omega_1 u_3 + \omega_2 v_3).$$

Also define

$$\Phi(u; v) := 2p_{(3)}(u; v) - p_{(2,1)}(u; v) + p_{(1,1,1)}(u; v),$$

which is precisely the left-hand side of (2.11). Then the nonnegativity of (2.11) is equivalent to the nonnegativity of $\Phi(u; v)$ under the restrictions that $u_i \geq v_i \geq 0$, $i = 1, 2, 3$. Equivalently, to establish (2.5), we only have to show that the polynomial $\Phi(u + v; v)$ is nonnegative under the restrictions $u_i \geq 0$, $v_i \geq 0$, $i = 1, 2, 3$. However, the package MAPLE produces the stronger result that, in the monomial expansion of $\Phi(u+v;v)$, all the coefficients appearing are nonnegative; in fact, MAPLE calculates the monomial expansion of $\Phi(u+v;v)$ to be

$$\begin{aligned}\Phi(u+v;v) &= \omega_1^2\omega_2(2v_1v_2v_3 + v_1v_2u_3 + v_1u_2v_3 + u_1v_2v_3) \\ &\quad + \omega_1\omega_2^2(v_1v_2v_3 + v_1v_2u_3 + v_1u_2v_3 + u_1v_2v_3).\end{aligned}$$

Now it becomes clear that generalizations of Theorem 1.1 can be established in a similar manner. Indeed, any conjectured generalization involving $m$ functions $f_1,\ldots,f_m$ will be valid if, with $u = (u_1,\ldots,u_m)$ and $v = (v_1,\ldots,v_m)$, the coefficients in the monomial expansion of the corresponding polynomial $\Phi(u+v;v)$ are all nonnegative. This is the approach we adopt to establish the nonnegativity of the fourth- and fifth-order conjugate cumulants.

PROOF OF THEOREM 1.3. To show that (1.8) is nonnegative, we follow the same approach as in the case of Theorem 1.1. To initiate the proof by induction, it is straightforward to verify that the result is valid for the case in which $A = \varnothing$.

Now we turn to the inductive hypothesis for nonempty $A$. In the case of four increasing functions $f_i$, $i = 1, 2, 3, 4$, the claim analogous to (2.5) is that, for any $B \subseteq A$,

$$\begin{aligned}(2.16) \quad & 6\mathbb{E}_B((f_1f_2f_3f_4)_B) - 2[\mathbb{E}_B((f_1f_2f_3)_Bf_{4B}) + \cdots] \\ & - [\mathbb{E}_B((f_1f_2)_B(f_3f_4)_B) + \cdots] \\ & + [\mathbb{E}_B((f_1f_2)_Bf_{3B}f_{4B}) + \cdots] + \mathbb{E}_B(f_{1B}f_{2B}f_{3B}f_{4B}) \geq 0.\end{aligned}$$

Proceeding as in Remark 2.4, we apply MAPLE to verify that all coefficients in the monomial expansion of the corresponding polynomial $\Phi(u+v;v)$ are nonnegative. Once this has been done, the remainder of the proof follows the arguments given in the latter part of the proof of Theorem 1.1.

To prove that (1.9) is nonnegative, we begin with the claim that

$$\begin{aligned}(2.17) \quad & 24\mathbb{E}_B((f_1f_2f_3f_4f_5)_B) - 2[\mathbb{E}_B((f_1f_2f_3f_4)_Bf_{5B}) + \cdots] \\ & - [\mathbb{E}_B((f_1f_2f_3)_B(f_4f_5)_B) + \cdots] \\ & - [\mathbb{E}_B((f_1f_2f_3)_Bf_{4B}f_{5B}) + \cdots] \\ & + [\mathbb{E}_B((f_1f_2)_B(f_3f_4)_Bf_{5B}) + \cdots] \\ & + [\mathbb{E}_B((f_1f_2)_Bf_{3B}f_{4B}f_{5B}) + \cdots] - \mathbb{E}_B(f_{1B}f_{2B}f_{3B}f_{4B}f_{5B}) \geq 0\end{aligned}$$

for any $B \subseteq A$. Next we construct the corresponding polynomial $\Phi(u;v)$, apply MAPLE to verify the nonnegativity of all coefficients in the monomial expansion of $\Phi(u+v;v)$, and then the remainder of the proof is as before. □

To formulate a conjecture for the case of $m$ increasing functions, we need some preliminaries from the theory of partitions [see Macdonald (1995)].

A *partition* $\lambda = (\lambda_1, \lambda_2, \ldots)$ is a sequence of nonnegative integers with $\lambda_1 \geq \lambda_2 \geq \cdots$. The *parts* of $\lambda$ are the nonzero $\lambda_i$; the *weight* of $\lambda$ is $|\lambda| := \lambda_1 + \lambda_2 + \cdots$; and the *length* of $\lambda$, denoted by $l(\lambda)$, is the number of parts of $\lambda$.

Given a partition $\lambda$, for each $i = 1, 2, \ldots$, let $\lambda_i'$ denote the cardinality of the set $\{j : \lambda_j \geq i\}$. Then $\lambda_1' \geq \lambda_2' \geq \cdots$, and the partition $\lambda' = (\lambda_1', \lambda_2', \ldots)$ is called the *partition conjugate to* $\lambda$. It is not difficult to verify that $(\lambda')' = \lambda$ and that $\lambda_1' = l(\lambda)$.

For $m \in \mathbb{N}$, $m \geq 2$, let $\mathfrak{S}_m$ denote the symmetric group on $m$ symbols. For any permutation $\tau \in \mathfrak{S}_m$ and any vector $u = (u_1, \ldots, u_m) \in \mathbb{R}^m$, define $\tau \cdot u := (u_{\tau(1)}, \ldots, u_{\tau(m)})$, the standard action of $\mathfrak{S}_m$ on $\mathbb{R}^m$.

Let $f_1, \ldots, f_m$ be functions on the lattice $L$ and let $\mu$ be a probability measure on $L$. For any partition $\lambda = (\lambda_1, \ldots, \lambda_m)$ of weight $m$, define

$$\mathcal{P}_\lambda(f_1, \ldots, f_m) := \prod_{j=1}^{l(\lambda)} \mathbb{E}\left(\prod_{k=1}^{\lambda_j} f_{\lambda_1 + \cdots + \lambda_{j-1} + k}\right),$$

where expectations are with respect to the measure $\mu$. We denote by $D(\lambda)$ the set of all $\tau \in \mathfrak{S}_m$ which give rise to *distinct* permutations $\mathcal{P}_\lambda(\tau \cdot (f_1, \ldots, f_m))$ of $\mathcal{P}_\lambda(f_1, \ldots, f_m)$. Then our conjecture for an $m$th-order generalization of the FKG inequality is the following:

CONJECTURE 2.5. *Let $\mu$ be an $MTP_2$ probability measure on the finite distributive lattice $L$, and let $f_1, \ldots, f_m$ be nonnegative increasing functions on $L$. Then*

(i) *There exists a set of nonzero constants $\{c_\lambda \in \mathbb{Z} : |\lambda| = m\}$ such that*

$$\mathcal{P}_m(f_1, \ldots, f_m) := \sum_{|\lambda|=m} c_\lambda \sum_{\tau \in D(\lambda)} \mathcal{P}_\lambda(\tau \cdot (f_1, \ldots, f_m))$$

*is nonnegative.*

(ii) *For each $m \geq 3$, there exists a constant $d_m$ such that*

(2.18) $$\mathcal{P}_m(f_1, \ldots, f_{m-1}, 1) \equiv d_m \mathcal{P}_{m-1}(f_1, \ldots, f_{m-1}).$$

(iii) *If $f_j \equiv 1$ for all $j = 1, \ldots, m$, then $\mathcal{P}_m(1, \ldots, 1) = 0$; equivalently,*

(2.19) $$\sum_{|\lambda|=m} \mathrm{card}(D(\lambda)) c_\lambda = 0.$$

For general $m$, the $m$th-order cumulant of a set of random variables $f_1, \ldots, f_m$ is

$$
\begin{aligned}
(2.20) \quad \kappa_m&(f_1, \ldots, f_m) \\
&:= \sum_{|\lambda|=m} (-1)^{l(\lambda)-1}(l(\lambda)-1)! \sum_{\tau \in D(\lambda)} \mathcal{P}_\lambda(\tau \cdot (f_1, \ldots, f_m)),
\end{aligned}
$$

and we define the $m$th-order conjugate cumulant

$$
\begin{aligned}
(2.21) \quad \kappa'_m&(f_1, \ldots, f_m) \\
&:= \sum_{|\lambda|=m} (-1)^{l(\lambda)-1}(l(\lambda')-1)! \sum_{\tau \in D(\lambda)} \mathcal{P}_\lambda(\tau \cdot (f_1, \ldots, f_m)).
\end{aligned}
$$

For $m = 2, 3, 4, 5$, the conjecture is valid if we choose for $\mathcal{P}_m$ the conjugate cumulant $\kappa'_m$. For these values of $m$, the coefficients $c_\lambda = (-1)^{l(\lambda)-1}(l(\lambda') - 1)!$, $|\lambda| = m$, in the expansion of $\kappa'_m(f_1, \ldots, f_m)$ are all nonzero and satisfy (2.19). Also, it can be verified that, for $m = 3, 4, 5$, the $\kappa'_m$ satisfy (2.18) with $d_m = m - 2$. In light of Theorems 1.1 and 1.3, all the preceding calculations provide evidence for the general conjecture.

For $m \geq 6$, it appears that the conjugate cumulants do not satisfy (2.19). In fact, it appears to us that $\kappa'_m(1, \ldots, 1) > 0$ for all $m \geq 6$.

**3. Diffusion processes.** Once we have a generalization of the FKG inequality within the context of finite distributive lattices, we can transfer that result to the context of functions and measures on $\mathbb{R}^n$ using standard approximation procedures [see Karlin and Rinott (1980)]. In this section we give a direct approach, using the method of diffusion processes, to our generalizations of the FKG inequalities on $\mathbb{R}^n$. This approach is based on the ideas of Herbst and Pitt (1991) and we first present some preliminary material, all of which is abstracted from Herbst and Pitt (1991).

We denote by $C(\mathbb{R}^n)$ the space of real-valued continuous functions on $\mathbb{R}^n$. Further, we define $C_b(\mathbb{R}^n) = \{f \in C(\mathbb{R}^n) : \|f\|_\infty < \infty\}$ and denote by $C_b^\infty(\mathbb{R}^n)$ the space of functions $f$ which have bounded continuous derivatives of all orders. We also denote by $\mathcal{M}$ the space of functions $f \in C_b(\mathbb{R}^n)$ which are increasing.

A probability measure $\mu$ on $\mathbb{R}^n$ is said to be *positively correlated*, or *associated*, if

$$\int_{\mathbb{R}^n} f_1(x) f_2(x) \, d\mu(x) \geq \int_{\mathbb{R}^n} f_1(x) \, d\mu(x) \cdot \int_{\mathbb{R}^n} f_2(x) \, d\mu(x)$$

for all $f_1, f_2 \in \mathcal{M}$. We denote by $\mathcal{P}$ the collection of all probability measures $\mu$ which are positively correlated.

Let $\{P_t : t \geq 0\} = \{P(t; x, dy) : t \geq 0\}$ be a Markov transition semigroup of probability measures on $\mathbb{R}^n$. Assume that $\{P_t : t \geq 0\}$ is Feller-continuous; that is, the operators $P_t$, defined by

$$P_t f(x) := \int_{\mathbb{R}^n} f(y) P(t; x, dy),$$

map the space $C_b(\mathbb{R}^n)$ into itself. The semigroup $\{P_t : t \geq 0\}$ is called *monotonic*, or is said to leave $\mathcal{M}$ *invariant*, if, for all $t \geq 0$, $P_t f \in \mathcal{M}$ whenever $f \in \mathcal{M}$. The semigroup $\{P_t : t \geq 0\}$ also acts on measures by

$$\mu P_t(A) := \int_{\mathbb{R}^n} \mu(dx) P(t; x, A) \equiv \int_{\mathbb{R}^n} \int_{\mathbb{R}^n} I_A(y) \mu(dx) P(t; x, dy),$$

where $I_A(\cdot)$ denotes the indicator function of the measurable set $A$; equivalently, for any $f \in C_b(\mathbb{R}^n)$, the measure $\mu P_t$ is defined by

$$\int_{\mathbb{R}^n} f(x) \, d\mu P_t(x) = \int_{\mathbb{R}^n} \int_{\mathbb{R}^n} f(y) \, d\mu(x) P(t; x, dy).$$

If $\mu P_t \in \mathcal{P}$ for all $t \geq 0$ whenever $\mu \in \mathcal{P}$, then we shall say that $\{P_t : t \geq 0\}$ *preserves positive correlations*. Thus $\{P_t : t \geq 0\}$ preserves positive correlations if

$$(3.1) \qquad \mu P_t(f_1 f_2) - (\mu P_t f_1)(\mu P_t f_2) \geq 0$$

for all $\mu \in \mathcal{P}$, $t \geq 0$ and $f_1, f_2 \in \mathcal{M}$.

The (*strong*) *infinitesimal generator* $\bar{G}$ of $\{P_t : t \geq 0\}$ is the linear operator on $C_b(\mathbb{R}^n)$ given by

$$(3.2) \qquad \bar{G} f(x) := \lim_{\varepsilon \to 0+} \frac{P_\varepsilon f(x) - f(x)}{\varepsilon},$$

where the convergence is uniform in $x$. The domain $\mathcal{D}(\bar{G})$ of $\bar{G}$ is the class of all functions $f$ for which the limit exists.

For each $i = 1, \ldots, n$, denote by $\partial_i$ the partial derivative, $\partial/\partial x_i$, with respect to $x_i$, the $i$th coordinate of the vector $x$.

Let $a(x) = (a_{i,j}(x))$, $x \in \mathbb{R}^n$, be a symmetric positive semidefinite matrix-valued function on $\mathbb{R}^n$ and let $b(x) = (b_j(x))$ be a vector field on $\mathbb{R}^n$. Denote by $G$ the differential operator

$$(3.3) \qquad G f(x) := \tfrac{1}{2} \sum_{i,j=1}^n a_{i,j}(x) \, \partial_i \, \partial_j f(x) + \sum_{j=1}^n b_j(x) \, \partial_j f(x)$$

with domain $\mathcal{D}(G) = \{f \in C_b^\infty(\mathbb{R}^n) : Gf \in C_b(\mathbb{R}^n)\}$.

Following Herbst and Pitt (1991), we call $\{P_t : t \geq 0\}$ a *diffusion semigroup with diffusion coefficients* $a(x)$ *and* $b(x)$ if $\mathcal{D}(G) \subseteq \mathcal{D}(\bar{G})$ and $Gf = \bar{G}f$ for all $f \in \mathcal{D}(G)$. We will assume throughout that $G$ *generates* $\{P_t : t \geq 0\}$ in the sense that there exists a unique diffusion semigroup $\{P_t : t \geq 0\}$ corresponding to given coefficients $a(x)$ and $b(x)$; a sufficient condition under which $G$ generates $\{P_t : t \geq 0\}$ is given by Herbst and Pitt [(1991), (1.3)] [cf. Chen and Wang (1993)].

Applying the semigroup property $P_{s+\varepsilon} = P_s P_\varepsilon$ and (3.2), it is straightforward to show that, for any $f \in \mathcal{D}(G)$,

$$(3.4) \quad \frac{\partial}{\partial s} P_s f(x) := \lim_{\varepsilon \to 0+} \frac{P_{s+\varepsilon} f(x) - P_s f(x)}{\varepsilon} = G P_s f(x) = P_s G f(x),$$

with the limit holding uniformly in $x$.

For $f_1, f_2 \in C^1(\mathbb{R}^n)$, we shall use the notation

$$(3.5) \qquad \Gamma_1(f_1, f_2)(x) := \sum_{i,j=1}^n a_{i,j}(x) \, \partial_i f_1(x) \, \partial_j f_2(x),$$

$x \in \mathbb{R}^n$; this operator is also known as the *carré du champ* operator [cf. Hu (2000)]. For $f_1, f_2 \in \mathcal{D}(G)$, it is straightforward to verify that

$$(3.6) \quad \Gamma_1(f_1, f_2)(x) = G(f_1 f_2)(x) - f_1(x) G(f_2)(x) - f_2(x) G(f_1)(x),$$

so that $\Gamma_1$ can be viewed as measuring the extent to which the operator $G$ is a derivation.

By combining the results of Herbst and Pitt (1991) and Chen and Wang (1993), we obtain the following criterion for monotonicity of the semigroup $\{P_t : t \geq 0\}$.

THEOREM 3.1 [Herbst and Pitt (1991) and Chen and Wang (1993)]. *The semigroup $\{P_t : t \geq 0\}$ is monotonic if and only if*

(i) *for all $i, j = 1, \ldots, n$, $a_{i,j}(x)$ depends only on $x_i$ and $x_j$, and*
(ii) *for all $j = 1, \ldots, n$, $b_j(x) \geq b_j(y)$ whenever $x \geq y$ with $x_j = y_j$.*

A necessary and sufficient condition for preservation of positive correlations is the following result of Herbst and Pitt [(1991), Theorem 1.3].

THEOREM 3.2 [Herbst and Pitt (1991)]. *The semigroup $\{P_t : t \geq 0\}$ preserves positive correlations if and only if*

(i) $\{P_t : t \geq 0\}$ *is monotonic, and*
(ii) $a_{i,j}(x) \geq 0$ *for all $x \in \mathbb{R}^n$, $i, j = 1, \ldots, n$.*

Theorem 3.1 can be viewed as providing conditions under which the semigroup $\{P_t : t \geq 0\}$ preserves "second-order" correlations. In the sequel we show, generalizing (3.1), how $\{P_t : t \geq 0\}$ preserves certain third-order correlations. The following result is motivated by Herbst and Pitt (1991), notably their Proposition 4.1, and the proof below is along the lines developed by them.

PROPOSITION 3.3. *Suppose that the diffusion semigroup $\{P_t : t \geq 0\}$ is monotonic, that*

$$\Gamma_1(f_1, f_2)(x) \geq 0 \tag{3.7}$$

*for all $x \in \mathbb{R}^n$ and for all smooth $f_1, f_2 \in \mathcal{M}$, and that $a_{i,j} \in \mathcal{M}$ for all $i, j = 1, \ldots, n$. Then*

$$\begin{aligned}
2P_t(f_1 f_2 f_3)(x) &- [P_t(f_1 f_2)(x) P_t f_3(x) \\
&+ P_t(f_1 f_3)(x) P_t f_2(x) + P_t f_1(x) P_t(f_2 f_3)(x)] \\
&+ P_t f_1(x) P_t f_2(x) P_t f_3(x) \geq 0
\end{aligned} \tag{3.8}$$

*for all nonnegative smooth $f_1, f_2, f_3 \in \mathcal{M}$ and all $x \in \mathbb{R}^n$.*

PROOF. Since the class of smooth functions is dense in $\mathcal{M}$ in the topology of bounded locally uniform convergence, it suffices to prove (3.8) for the case in which $f_1, f_2, f_3 \in \mathcal{M} \cap C_b^\infty(\mathbb{R}^n)$.

Now fix $t > 0$ and define the function

$$\begin{aligned}
h(s) = {}& 2P_s(P_{t-s}f_1 \cdot P_{t-s}f_2 \cdot P_{t-s}f_3)(x) \\
& - [P_s(P_{t-s}f_1 \cdot P_{t-s}f_2)(x) \cdot P_t f_3(x) \\
& \quad + P_s(P_{t-s}f_1 \cdot P_{t-s}f_3)(x) \cdot P_t f_2(x) \\
& \quad + P_t f_1(x) \cdot P_s(P_{t-s}f_2 \cdot P_{t-s}f_3)(x)] \\
& + P_t f_1(x) \cdot P_t f_2(x) \cdot P_t f_3(x),
\end{aligned} \tag{3.9}$$

$0 \leq s \leq t$. Observe that $h(0) = 0$ and

$$\begin{aligned}
h(t) = {}& 2P_t(f_1 f_2 f_3)(x) \\
& - [P_t(f_1 f_2)(x) \cdot P_t f_3(x) \\
& \quad + P_t(f_1 f_3)(x) \cdot P_t f_2(x) + P_t f_1(x) \cdot P_t(f_2 f_3)(x)] \\
& + P_t f_1(x) \cdot P_t f_2(x) \cdot P_t f_3(x).
\end{aligned} \tag{3.10}$$

Thus, to establish (3.8), we need only show that $h'(s) \geq 0$ for all $s \in (0, t)$.

We now apply (3.4) repeatedly to (3.9) to differentiate $h$. Suppressing the notational dependence of all functions on the argument $x$, we obtain

$$h'(s) = 2[P_s G(P_{t-s}f_1 \cdot P_{t-s}f_2 \cdot P_{t-s}f_3) - P_s(GP_{t-s}f_1 \cdot P_{t-s}f_2 \cdot P_{t-s}f_3)$$
$$- P_s(P_{t-s}f_1 \cdot GP_{t-s}f_2 \cdot P_{t-s}f_3) - P_s(P_{t-s}f_1 \cdot P_{t-s}f_2 \cdot GP_{t-s}f_3)]$$
$$- [P_s G(P_{t-s}f_1 \cdot P_{t-s}f_2)$$
$$- P_s(GP_{t-s}f_1 \cdot P_{t-s}f_2) - P_s(P_{t-s}f_1 \cdot GP_{t-s}f_2)]P_t f_3$$
$$- [P_s G(P_{t-s}f_1 \cdot P_{t-s}f_3)$$
$$- P_s(GP_{t-s}f_1 \cdot P_{t-s}f_3) - P_s(P_{t-s}f_1 \cdot GP_{t-s}f_3)]P_t f_2$$
$$- [P_s G(P_{t-s}f_2 \cdot P_{t-s}f_3)$$
$$- P_s(GP_{t-s}f_2 \cdot P_{t-s}f_3) - P_s(P_{t-s}f_2 \cdot GP_{t-s}f_3)]P_t f_1.$$

Equivalently, denoting $P_{t-s}f_i$ by $g_i$, $i=1,2,3$, what we have shown is that

$$h'(s) = 2[P_s G(g_1 g_2 g_3) - P_s(g_2 g_3 G g_1) - P_s(g_1 g_3 G g_2) - P_s(g_1 g_2 G g_3)]$$
$$- [P_s G(g_1 g_2) - P_s(g_2 G g_1) - P_s(g_1 G g_2)]P_s g_3$$
$$- [P_s G(g_1 g_3) - P_s(g_3 G g_1) - P_s(g_1 G g_3)]P_s g_2$$
$$- [P_s G(g_2 g_3) - P_s(g_3 G g_2) - P_s(g_2 G g_3)]P_s g_1.$$

Using the definition of $\Gamma_1(\cdot,\cdot)$ in (3.5), we have

$$h'(s) = 2P_s[G(g_1 g_2 g_3) - g_2 g_3 G g_1 - g_1 g_3 G g_2 - g_1 g_2 G g_3]$$
(3.11)
$$- [P_s \Gamma_1(g_1, g_2) \cdot P_s g_3$$
$$+ P_s \Gamma_1(g_1, g_3) \cdot P_s g_2 + P_s \Gamma_1(g_2, g_3) \cdot P_s g_1].$$

For $\{i,j,k\} = \{1,2,3\}$, we express each term $P_s \Gamma_1(g_i, g_j) \cdot P_s g_k$ in (3.11) in the form

$$P_s \Gamma_1(g_i, g_j) \cdot P_s g_k$$
$$\equiv P_s(\Gamma_1(g_i, g_j) \cdot g_k) - [P_s(\Gamma_1(g_i, g_j) \cdot g_k) - P_s \Gamma_1(g_i, g_j) \cdot P_s g_k],$$

and then we obtain

$$h'(s) = P_s \Gamma_1(g_1, g_2, g_3) + [P_s(\Gamma_1(g_1, g_2) \cdot g_3) - P_s \Gamma_1(g_1, g_2) \cdot P_s g_3]$$
(3.12)
$$+ [P_s(\Gamma_1(g_1, g_3) \cdot g_2) - P_s \Gamma_1(g_1, g_3) \cdot P_s g_2]$$
$$+ [P_s(\Gamma_1(g_2, g_3) \cdot g_1) - P_s \Gamma_1(g_2, g_3) \cdot P_s g_1],$$

where

$$\Gamma_1(g_1, g_2, g_3) := 2[G(g_1 g_2 g_3) - g_2 g_3 G g_1 - g_1 g_3 G g_2 - g_1 g_2 G g_3]$$
(3.13)
$$- [g_3 \Gamma_1(g_1, g_2) + g_2 \Gamma_1(g_1, g_3) + g_1 \Gamma_1(g_2, g_3)]$$
$$= 2G(g_1 g_2 g_3) - g_1 G(g_2 g_3) - g_2 G(g_1 g_3) - g_3 G(g_1 g_2).$$

By a direct calculation using (3.3) and (3.5), we obtain from (3.13) the identity

$$\Gamma_1(g_1, g_2, g_3) = g_1 \sum_{i,j=1}^{n} a_{i,j}\, \partial_i g_2\, \partial_j g_3$$

(3.14)
$$+ g_2 \sum_{i,j=1}^{n} a_{i,j}\, \partial_i g_1\, \partial_j g_3 + g_3 \sum_{i,j=1}^{n} a_{i,j}\, \partial_i g_1\, \partial_j g_2$$

$$\equiv g_1 \Gamma_1(g_2, g_3) + g_2 \Gamma_1(g_1, g_3) + g_3 \Gamma_1(g_1, g_2).$$

Since $f_1$, $f_2$ and $f_3$ are nonnegative, then so are $g_1$, $g_2$ and $g_3$. By the assumptions that $\{P_t : t \geq 0\}$ is monotonic and that (3.7) holds, we find that $g_i \in \mathcal{M} \cap C_b^\infty(\mathbb{R}^n)$, $i = 1, 2, 3$, and $\Gamma_1(g_k, g_l) \geq 0$, $1 \leq k < l \leq 3$; therefore $\Gamma_1(g_1, g_2, g_3) \geq 0$.

Since each $g_k \in \mathcal{M}$, then $\partial_i g_k \geq 0$ for all $i = 1, \ldots, n$. By assumption, each $a_{i,j} \in \mathcal{M}$. Therefore it follows from (3.5) that, for any pair of smooth $g_k, g_l \in \mathcal{M}$, $\Gamma_1(g_k, g_l)$ is a nonnegative linear combination of elements of $\mathcal{M}$; hence $\Gamma_1(g_i, g_j) \in \mathcal{M}$. By Proposition 4.1 of Herbst and Pitt (1991), the semigroup $\{P_t : t \geq 0\}$ is known to preserve positive correlations; therefore each term $P_s(\Gamma_1(g_i, g_j) \cdot g_k) - P_s \Gamma_1(g_i, g_j) \cdot P_s g_k$ is nonnegative, establishing that (3.12) is a decomposition of $h'(s)$ into nonnegative terms. This proves that $h'(s) \geq 0$, $s \in (0, t)$, from which we conclude that $h(t) \geq 0$. □

REMARK 3.4. (i) Consider the case in which $f_3 \equiv 1$ in the foregoing proof. Then the terms $P_s(\Gamma_1(g_i, g_j) \cdot g_k) - P_s \Gamma_1(g_i, g_j) \cdot P_s g_k$ are identically zero, so that no monotonicity conditions on the $a_{i,j}$ are required. Then we recover a result of Herbst and Pitt [(1991), Proposition 4.1].

(ii) By analogous arguments we can extend some results of Wang and Yan (1994) on positive correlations for diffusion processes on $n$-dimensional tori.

(iii) Results similar to Proposition 3.3 have played a prominent role in the study of functional inequalities, for example, log-Sobolev inequalities; see Hu (2000). In work now in progress, we will study extensions of those inequalities by means of generalizations of Proposition 3.3.

As a consequence of Proposition 3.3, we now obtain an analog of Theorem 1.1 for functions on $\mathbb{R}^n$; here again, we follow arguments given by Herbst and Pitt [(1991), Corollary 1.7] in their proof of the FKG inequality.

COROLLARY 3.5. *Suppose that $\mu(dx) = \exp(\psi(x))\, dx$ is a probability measure on $\mathbb{R}^n$, where $\psi \in C^2(\mathbb{R}^n)$ satisfies the properties that $\psi(x) > -\infty$ and $\partial_i \partial_j \psi(x) \geq 0$ for all $i \neq j$ and for all $x \in \mathbb{R}^n$. If $f_1$, $f_2$ and $f_3$ are nonnegative functions in $\mathcal{M}$ such that the various expectations exist, then $\kappa'_3(f_1, f_2, f_3) \geq 0$.*

PROOF. Construct the diffusion process $\{P_t : t \geq 0\}$ with infinitesimal generator

$$G = \sum_{j=1}^n \partial_j^2 + \sum_{j=1}^n (\partial_j \psi)\, \partial_j,$$

corresponding to (3.3) with $a_{i,j} = 2\delta_{i,j}$ and $b_j = \partial_j \psi$. It is well known that the process $\{P_t : t \geq 0\}$ is monotonic [by Theorem 3.1, this is equivalent to the assumption that $\partial_i \partial_j \psi(x) \geq 0$ for all $x$ and $i \neq j$]. Also, $\Gamma_1(f_1, f_2) \geq 0$ for all smooth $f_1, f_2 \in \mathcal{M}$ and the $a_{i,j}$, being constants, clearly are elements of $\mathcal{M}$. Applying Proposition 3.3, we deduce that (3.8) holds. Letting $t \to \infty$ in (3.8) and noting that $P_t f \to \mathbb{E}(f)$, almost surely, we obtain (1.7). $\square$

It also is clear how to extend Proposition 3.3 to the case of four or more functions. We have found the algebraic calculations to be more extensive, indicating that an alternative approach is needed to resolve the case of an arbitrary number of functions. Nevertheless, the strategy is the same, and we illustrate it by sketching the details for the case of four functions. Suppose that $f_1, \ldots, f_4 \in \mathcal{M}$ are nonnegative and smooth. For $x \in \mathbb{R}^n$, define $g_i = P_{t-s} f_i$, $i = 1, \ldots, 4$ and

$$\begin{aligned}
h(s) = {} & 6 P_s(g_1 g_2 g_3 g_4) \\
& - 2[P_s(g_1 g_2 g_3) \cdot P_t f_4 + P_s(g_1 g_2 g_4) \cdot P_t f_3 \\
& \qquad + P_s(g_1 g_3 g_4) \cdot P_t f_2 + P_s(g_2 g_3 g_4) \cdot P_t f_1] \\
& - [P_s(g_1 g_2) \cdot P_s(g_3 g_4) + P_s(g_1 g_3) \cdot P_s(g_2 g_4) + P_s(g_1 g_4) \cdot P_s(g_2 g_3)] \\
& + [P_s(g_1 g_2) \cdot P_t f_3 \cdot P_t f_4 + P_s(g_1 g_3) \cdot P_t f_2 \cdot P_t f_4 \\
& \qquad + P_s(g_1 g_4) \cdot P_t f_2 \cdot P_t f_3 + P_s(g_2 g_3) \cdot P_t f_1 \cdot P_t f_4 \\
& \qquad + P_s(g_2 g_4) \cdot P_t f_1 \cdot P_t f_3 + P_s(g_3 g_4) \cdot P_t f_1 \cdot P_t f_2] \\
& + P_t f_1 \cdot P_t f_2 \cdot P_t f_3 \cdot P_t f_4,
\end{aligned}$$

$0 \leq s \leq t$, where, as usual, we have suppressed all notational dependence of functions on $x$. Define

$$\begin{aligned}
& \Gamma_1(g_1, g_2, g_3, g_4) \\
& \quad := g_1 g_2 \Gamma_1(g_3, g_4) + g_1 g_3 \Gamma_1(g_2, g_4) + g_1 g_4 \Gamma_1(g_2, g_3) \\
& \qquad + g_2 g_3 \Gamma_1(g_1, g_4) + g_2 g_4 \Gamma_1(g_1, g_3) + g_3 g_4 \Gamma_1(g_1, g_2),
\end{aligned}$$

and denote $\Gamma_1(g_2, g_3, g_4)$ by $\Gamma_1(\{g_1, \ldots, g_4\} \setminus \{g_1\})$, and so on. For $0 < s < t$,

a lengthy, but straightforward, calculation leads to the result

$$h'(s) = P_s \Gamma_1(g_1, g_2, g_3, g_4)$$
$$+ 2 \sum_{j=1}^{4} [P_s(g_j \cdot \Gamma_1(\{g_1, \ldots, g_4\} \setminus \{g_j\}))$$
$$- P_s g_j \cdot P_s \Gamma_1(\{g_1, \ldots, g_4\} \setminus \{g_j\})]$$
$$+ \sum_{1 \leq i < j \leq 4} [P_s(g_i g_j \cdot \Gamma_1(\{g_1, \ldots, g_4\} \setminus \{g_i, g_j\}))$$
$$- P_s(g_i g_j) \cdot P_s \Gamma_1(\{g_1, \ldots, g_4\} \setminus \{g_i, g_j\})]$$
$$+ \sum_{1 \leq i < j \leq 4} P_s g_i \cdot P_s g_j \cdot P_s \Gamma_1(\{g_1, \ldots, g_4\} \setminus \{g_i, g_j\}).$$

Now we apply the same arguments which concluded the proof of Proposition 3.3. Then the above expression is seen to be a resolution of $h'(s)$ into nonnegative terms; hence, under the hypotheses of Proposition 3.3, $h'(s) \geq 0$, $0 < s < t$, and therefore $h(t) \geq h(0) = 0$. Finally, we argue as in Corollary 3.5 to deduce $\kappa_4' \geq 0$.

It is interesting to compare the techniques applied in Sections 2 and 3. On the one hand, the lattice formulation in Section 2 required no analytical considerations; however, the algebraic manipulations were sufficiently complicated that we found it necessary to use the MAPLE package to decompose expressions as sums of nonnegative terms. On the other hand, in the present section, analytical machinery was needed but the algebraic calculations appeared to be simpler and we made no use of the MAPLE package.

**4. Duplicate variables.** Another approach to establishing the classical FKG inequality on $\mathbb{R}^n$ is by way of the method of duplicate variables; see Cartier (1974) and Glimm and Jaffe (1987). Therefore it is natural to search for a proof of Theorem 1.1 using this method. For $n = 1$, we have such a proof; however, the inductive step seems difficult and we have not been able to find it. To demonstrate the difficulties inherent in the inductive step, we now establish the case in which $n = 1$.

Here, the functions $f_1$, $f_2$ and $f_3$ are nonnegative and increasing on $\mathbb{R}$ and we wish to show that the conjugate cumulant $\kappa_3'(f_1, f_2, f_3)$ is nonnegative. By duplicating variables we have

$$\kappa_3'(f_1, f_2, f_3) = \int_{\mathbb{R}^3} I(x_1, x_2, x_3) \, d\mu(x_1) \, d\mu(x_2) \, d\mu(x_3),$$

where

$$I(x_1, x_2, x_3)$$
$$= 2f_1(x_1)f_2(x_1)f_3(x_1)$$
$$- [f_1(x_1)f_2(x_1)f_3(x_2) + f_1(x_1)f_2(x_2)f_3(x_1) + f_1(x_2)f_2(x_1)f_3(x_1)]$$
$$+ f_1(x_1)f_2(x_2)f_3(x_3).$$

Define

$$J(x_1, x_2, x_3) = \sum_{\tau \in \mathfrak{S}_3} I(\tau \cdot (x_1, x_2, x_3));$$

then, by a symmetry argument, we obtain

$$\kappa'_3(f_1, f_2, f_3) = \frac{1}{3!} \int_{\mathbb{R}^3} J(x_1, x_2, x_3) \, d\mu(x_1) \, d\mu(x_2) \, d\mu(x_3).$$

With the help of MAPLE we find that

$$\begin{aligned}J(x_1, x_2, x_3) = {} & f_3(x_1)(f_1(x_3) - f_1(x_1))(f_2(x_3) - f_2(x_1)) \\
& + f_1(x_3)(f_2(x_3) - f_2(x_1))(f_3(x_3) - f_3(x_1)) \\
& + f_2(x_1)(f_1(x_3) - f_1(x_1))(f_3(x_3) - f_3(x_1)) \\
& + f_3(x_1)(f_1(x_2) - f_1(x_1))(f_2(x_2) - f_2(x_1)) \\
& + f_1(x_2)(f_2(x_2) - f_2(x_1))(f_3(x_2) - f_3(x_1)) \\
& + f_2(x_1)(f_1(x_2) - f_1(x_1))(f_3(x_2) - f_3(x_1)) \\
& + f_1(x_3)(f_2(x_3) - f_2(x_2))(f_3(x_3) - f_3(x_2)) \\
& + f_3(x_2)(f_1(x_3) - f_1(x_2))(f_2(x_3) - f_2(x_2)) \\
& + f_2(x_3)(f_1(x_3) - f_1(x_2))(f_3(x_3) - f_3(x_2)) \\
& + (f_1(x_3) - f_1(x_2))(f_2(x_2) - f_2(x_1))(f_3(x_3) - f_3(x_2)) \\
& + (f_1(x_3) - f_1(x_1))(f_2(x_3) - f_2(x_2))(f_3(x_3) - f_3(x_2)) \\
& + (f_1(x_3) - f_1(x_2))(f_2(x_3) - f_2(x_2))(f_3(x_2) - f_3(x_1)). \end{aligned} \quad (4.1)$$

Since $f_1$, $f_2$ and $f_3$ are nonnegative and increasing, we deduce immediately that $J(x_1, x_2, x_3) \geq 0$ on the "fundamental" chamber $\{x_1 < x_2 < x_3\}$. It then follows by symmetry that $J(x_1, x_2, x_3) \geq 0$ on every Weyl chamber $\{x_{\tau(1)} < x_{\tau(2)} < x_{\tau(3)}\}$, $\tau \in \mathfrak{S}_3$. By taking limits as $(x_1, x_2, x_3)$ goes to a boundary $\{x_{\tau(1)} \leq x_{\tau(2)} \leq x_{\tau(3)}\}$, $\tau \in \mathfrak{S}_3$, we find that $J(x_1, x_2, x_3)$ remains nonnegative on every wall of a chamber. Therefore $\kappa'_3(f_1, f_2, f_3) \geq 0$.

To carry out the inductive step toward higher dimensions requires a method for handling functions of the form $J$ in (4.1); unfortunately, we have not been able to develop such a technique. More importantly, our usage of the terminology "Weyl chamber," "wall" and "fundamental chamber" is not accidental, for we believe that inequalities of FKG-type exist within the context of arbitrary finite reflection groups; see Gross and Richards (1995).

**5. Applications.** Any generalization of the FKG inequality has the potential for vast numbers of applications. In this section we provide a few applications and point the reader to others. We begin with generalizations of some applications given by Seymour and Welsh (1975).

5.1. *Generalized inequalities for Bernstein polynomials.* For $f \in C[0,1]$, the Bernstein polynomial of $f$ is the polynomial defined on $[0,1]$ by

$$B_n f(x) = \sum_{k=0}^{n} f(k/n) \binom{n}{k} x^k (1-x)^{n-k}.$$

Generalizing Theorem 2.6 of Seymour and Welsh (1975), we have the following result.

PROPOSITION 5.1. *If $f_1$, $f_2$ and $f_3$ are nonnegative increasing functions on $[0,1]$, then their Bernstein polynomials satisfy*

$$\begin{aligned}
(5.1) \quad & 2B_n(f_1 f_2 f_3)(x) \\
& - [B_n(f_1 f_2)(x) B_n f_3(x) \\
& \quad + B_n(f_1 f_3)(x) B_n f_2(x) + B_n f_1(x) B_n(f_2 f_3)(x)] \\
& + B_n f_1(x) B_n f_2(x) B_n f_3(x) \geq 0
\end{aligned}$$

*for all $x \in [0,1]$.*

PROOF. Let $A$ be a set of $n$ elements and let $x \in [0,1]$. Define a probability measure $\mu$ on the lattice $2^A$ by

$$\mu(a) = c x^{\mathrm{card}(a)} (1-x)^{n-\mathrm{card}(a)},$$

$a \in 2^A$, where $\mathrm{card}(a)$ denotes the cardinality of $a$ and $c$ is the normalizing constant; it is well known that $\mu$ satisfies the $\mathrm{MTP}_2$ condition (2.1). Next we define

$$\widehat{f}_j(a) = f_j(\mathrm{card}(a)/n),$$

$j=1,2,3$, $a \in 2^A$. Then the functions $\widehat{f}_j$ are nonnegative and increasing on $2^A$. On applying Theorem 1.1 to the probability measure $\mu$ and the functions $\widehat{f}_j$, $j=1,2,3$, we obtain (5.1). □

By arguing along similar lines, inequalities for Bernstein polynomials for four and five functions can be obtained through application of Theorem 1.3.

5.2. *Inequalities for log-convex sequences.* Recall that a sequence of real numbers $\{a_0, \ldots, a_n\}$ is called *log-convex* if $a_k^2 \leq a_{k-1}a_{k+1}$ for all $k = 1, \ldots, n-1$.

For any sequence $\{\alpha_0, \ldots, \alpha_n\}$, we will use the notation

$$\langle \{\alpha_k\} \rangle := \sum_{k=0}^{n} \alpha_k.$$

The following result extends Theorem 3.2 of Seymour and Welsh (1975), which itself is a generalization of a classical inequality of Tchebycheff.

PROPOSITION 5.2. *Suppose that $\{a_0, \ldots, a_n\}$ is a positive, log-convex sequence with $\langle \{a_k\} \rangle = 1$. Suppose also that the sequences $\{\alpha_0, \ldots, \alpha_n\}$, $\{\beta_0, \ldots, \beta_n\}$ and $\{\gamma_0, \ldots, \gamma_n\}$ are increasing and nonnegative. Then*

$$\begin{aligned}
& 2\langle \{a_k \alpha_k \beta_k \gamma_k\} \rangle - \langle \{a_k \alpha_k \beta_k\} \rangle \langle \{a_k \gamma_k\} \rangle \\
& \quad - \langle \{a_k \alpha_k \gamma_k\} \rangle \langle \{a_k \beta_k\} \rangle - \langle \{a_k \alpha_k\} \rangle \langle \{a_k \beta_k \gamma_k\} \rangle \\
& \quad + \langle \{a_k \alpha_k\} \rangle \langle \{a_k \beta_k\} \rangle \langle \{a_k \gamma_k\} \rangle \geq 0.
\end{aligned} \tag{5.2}$$

PROOF. The proof of this result is similar to the proof of the corresponding result given by Seymour and Welsh (1975). Define the sequence

$$b_k = a_k \Big/ \binom{n}{k},$$

$k = 0, \ldots, n$. Since the sequences $\{a_k\}$ and $\{1/\binom{n}{k}\}$ are log-convex, then so is $\{b_k\}$.

Let $A = \{1, \ldots, n\}$ and define $\mu: 2^A \to \mathbb{R}$ by

$$\mu(a) = b_{\mathrm{card}(a)},$$

$a \in 2^A$. It is shown by Seymour and Welsh (1975), that $\mu$ satisfies (2.1).

Define $f_1, f_2, f_3: 2^A \to \mathbb{R}$ by

$$f_1(a) = \alpha_{\mathrm{card}(a)}, \qquad f_2(a) = \beta_{\mathrm{card}(a)}, \qquad f_3(a) = \gamma_{\mathrm{card}(a)},$$

$a \in 2^A$. Then the functions $f_j$ are nonnegative and increasing on $2^A$. On applying Theorem 1.1 to the probability measure proportional to $\mu$ and the functions $f_j$, $j = 1, 2, 3$, we obtain the desired result. $\square$

5.3. *A generalization of Kleitman's lemma.* A collection $C$ of subsets of a set $A$ is *closed above* if $a \in C$ and $a \subseteq b$ imply $b \in C$. Similarly, $C$ is *closed below* if $a \in C$ and $a \supseteq b$ imply $b \in C$.

Suppose that $A$ is a finite set of cardinality $n$, and let $U$ and $L$ be collections of subsets of $A$ such that $U$ is closed above and $L$ is closed below. Then Kleitman (1966) proved the remarkable inequality,

$$2^n \mathrm{card}(U \cap L) \leq \mathrm{card}(U) \mathrm{card}(L). \tag{5.3}$$

HIGHER-ORDER INEQUALITIES 25The following result generalizes this inequality; see Seymour and Welsh [(1975), Theorem 4.2].

PROPOSITION 5.3. *Let $A$ be a finite set of cardinality $n$, and let $U_1$, $U_2$, and $L$ be collections of subsets of $A$. Suppose that $U_1$ and $U_2$ are closed above and $L$ is closed below. Then*

(5.4)
$$2^{2n} \operatorname{card}(U_1 \cap U_2) - 2^{2n+1} \operatorname{card}(U_1 \cap U_2 \cap L)$$
$$+ 2^n ( \operatorname{card}(U_1 \cap U_2) \operatorname{card}(L)$$
$$+ \operatorname{card}(U_1 \cap L) \operatorname{card}(U_2) + \operatorname{card}(U_1) \operatorname{card}(U_2 \cap L))$$
$$- 2^n \operatorname{card}(U_1) \operatorname{card}(U_2) - \operatorname{card}(U_1) \operatorname{card}(U_2) \operatorname{card}(L) \geq 0.$$

PROOF. Let $\mu$ be the uniform distribution on $2^A$, that is, $\mu(a) = 1/2^n$ for all $a \in 2^A$. Let $f_1$ and $f_2$ be the characteristic functions of the sets $U_1$ and $U_2$, respectively; that is, for $j = 1, 2$,

$$f_j(a) = \begin{cases} 1, & \text{if } a \in U_j, \\ 0, & \text{if } a \notin U_j. \end{cases}$$

Also let $f_3$ denote the characteristic function of $L^c$, the complement of $L$. Then $f_1$, $f_2$ and $1 - f_3$ are nonnegative increasing functions. Then (5.4) is obtained by applying Theorem 1.1 to the measure $\mu$ and the functions $f_1$, $f_2$ and $1 - f_3$. □

For the case in which $U_1 \equiv U$ is closed above and $U_2 = 2^A$, the inequality (5.4) reduces to (5.3).

5.4. *Inequalities for matrix functions.*

DEFINITION 5.4. Let $R = (R(i,j))$ be an $n \times n$ real matrix. Then

(i) $R$ satisfies the *triangle property* if

(5.5) $$R(i,j)R(k,k) = R(i,k)R(k,j)$$

for all $i \leq k \leq j$ and all $j \leq k \leq i$ [see Barrett and Feinsilver (1981)];

(ii) $R$ is *nonnegative* if $R(i,j)$ is nonnegative for all $i$ and $j$;

(iii) $R$ *represents a discrete probability distribution* if $R$ is nonnegative and $\sum_{i,j=1}^n R(i,j) = 1$;

(iv) $R$ is *increasing* if $R(i,j)$ is monotone increasing in $i$ and in $j$.

For any $n \times n$ matrix $R = R(i,j)$ representing a discrete probability distribution and any $n \times n$ matrix $F = (F(i,j))$, we denote $\operatorname{tr}(RF')$ by $\mathbb{E}_R(F)$. Because

$$\operatorname{tr}(RF') = \sum_{i=1}^n \sum_{j=1}^n R(i,j) F(i,j),$$

the notation $\mathbb{E}_R(F)$ for $\mathrm{tr}(RF')$ is consistent with the notation for expectation with respect to the discrete probability distribution represented by the entries of the matrix $R$.

Given two $n \times n$ matrices $F_1 = (f_1(i,j))$ and $F_2 = (f_2(i,j))$, recall that the *Hadamard product* of $F_1$ and $F_2$ is the matrix $F_1 \circ F_2 = (f_1(i,j)f_2(i,j))$.

Now we have the following result.

PROPOSITION 5.5. *Let $R$ be an $n \times n$ matrix which represents a discrete probability distribution and satisfies the triangle property. If $F_1 = (f_1(i,j))$, $F_2 = (f_2(i,j))$ and $F_3 = (f_3(i,j))$ are nonnegative increasing $n \times n$ matrices, then*

$$
\begin{aligned}
2\mathbb{E}_R(F_1 \circ F_2 \circ F_3) &- [\mathbb{E}_R(F_1 \circ F_2)\mathbb{E}_R(F_3) \\
&\quad + \mathbb{E}_R(F_1 \circ F_3)\mathbb{E}_R(F_2) + \mathbb{E}_R(F_1)\mathbb{E}_R(F_2 \circ F_3)] \\
&+ \mathbb{E}_R(F_1)\mathbb{E}_R(F_2)\mathbb{E}_R(F_3) \geq 0.
\end{aligned}
\tag{5.6}
$$

*In particular,*

$$
\mathbb{E}_R(F_1 \circ F_2) \geq \mathbb{E}_R(F_1)\mathbb{E}_R(F_2).
\tag{5.7}
$$

PROOF. Let $A = \{1,\ldots,n\}$, so that $A \times A$ is the set of pairs of positive integers ranging from 1 through $n$. Equip $A \times A$ with the partial order given by $(i,j) \preccurlyeq (k,l)$ if $i \leq j$ and $k \leq l$; this is the same partial ordering utilized in the study of the FKG inequality on Euclidean space. For pairs $p = (i,j)$ and $q = (k,l)$ in $A \times A$, we define $p \vee q$ and $p \wedge q$ in the usual component-wise manner:

$$p \vee q = (\max(i,k), \max(j,l)), \qquad p \wedge q = (\min(i,k), \min(j,l)).$$

It is straightforward to check that the triangle property (5.5) is equivalent to the $\mathrm{MTP}_2$ condition with equality:

$$R(p \vee q)R(p \wedge q) = R(p)R(q)$$

for all $p, q \in A \times A$. Therefore, any matrix $R$ which represents a discrete distribution and satisfies the triangle condition corresponds to an $\mathrm{MTP}_2$ probability distribution on $A \times A$. This observation is the crux of the proof, for we can now apply Theorem 1.1 to nonnegative increasing matrix functions to deduce (5.6). Finally, (5.7) is the special case of (5.6) in which $f_3(i,j) \equiv 1$. □

We now turn to inequalities for rank and determinant functions of positive-definite functions.

Let $M$ denote a fixed $n \times n$ positive-semidefinite matrix $M$. Given an index set $a \subseteq A = \{1, \ldots, n\}$, we denote by $M[a]$ the submatrix of $M$ appearing in the rows and columns labelled by the elements of the set $a$. For any $a, b \subseteq A$, it is a result of Lundquist and Barrett (1996) that

(5.8) $\quad \operatorname{rank} M[a \cup b] + \operatorname{rank} M[a \cap b] \leq \operatorname{rank} M[a] + \operatorname{rank} M[b].$

For fixed $t > 0$, define the probability measure $\mu$ on $2^A$ by

$$\mu(a) = \frac{t^{n-\operatorname{rank} M[a]}}{\sum_{\alpha \subseteq A} t^{n-\operatorname{rank} M[\alpha]}},$$

$a \in 2^A$; the measure $\mu$ can be viewed as a generating function for the ranks of the submatrices $M[a]$, $a \subseteq A$.

It follows from (5.8) that $\mu$ satisfies the $\mathrm{MTP}_2$ condition (1.1) on the lattice $2^A$. By applying Theorem 1.1 to the measure $\mu$ and nonnegative increasing functions $f_j$, $j = 1, 2, 3$, we obtain various positivity results. We leave it to the reader to work out special cases, for example, the case in which the $f_j$ are defined in terms of the characteristic functions of subsets of $2^A$ that are closed above or below, as necessary.

Other examples arise by specifying the $f_j$ to be functions analogous to those chosen by Seymour and Welsh (1975) in the proof of their Theorem 5.10. More generally, we may consider the context of matroid theory considered by Seymour and Welsh (1975) and deduce by application of Theorem 1.1 higher-order total positivity properties of the rank-generating function of a matroid.

To derive determinantal inequalities, suppose that $M$ is a fixed positive-definite symmetric $n \times n$ matrix. We recall the generalized Hadamard–Fischer inequality

(5.9) $\quad \det M[a \cup b] \det M[a \cap b] \leq \det M[a] \det M[b];$

see Horn and Johnson [(1985), page 485]. For fixed $t > 0$, define the probability measure

$$\mu(a) = \frac{\det M[a]^{-t}}{\sum_{\alpha \subseteq A} \det M[\alpha]^{-t}},$$

$a \in 2^A$. Then $\mu$ satisfies the $\mathrm{MTP}_2$ condition (1.1) on the lattice $2^A$. On applying Theorem 1.1 to various choices of the functions $f_j$, we obtain determinantal inequalities.

For example, suppose $f_1(a) = \operatorname{tr} M[a]$, the trace of $M[a]$; $f_2(a) = \lambda_{\min}(M[a])$, the smallest eigenvalue of $M[a]$; and $f_3(a) = 1/\lambda_{\max}(M[a])$, the inverse of the largest eigenvalue of $M[a]$. Then it is a consequence of the fundamental *inclusion principle* describing the interlacing properties of eigenvalues of submatrices of Hermitian matrices [Horn and Johnson (1985), Theorem 4.3.15,

page 185] that these $f_j$ are all increasing functions on the lattice $2^A$. Therefore we may obtain eigenvalue inequalities by application of Theorem 1.1. As a special case, by applying the FKG inequality to the functions $f_2$ and $f_3$, we obtain the inequality

$$\left(\sum_{a\subseteq A} \det M[a]^{-t}\right)\left(\sum_{a\subseteq A} \det M[a]^{-t}\frac{\lambda_{\min}(M[a])}{\lambda_{\max}(M[a])}\right)$$
$$\geq \left(\sum_{a\subseteq A} \det M[a]^{-t}\lambda_{\min}(M[a])\right)\left(\sum_{a\subseteq A} \det M[a]^{-t}\frac{1}{\lambda_{\max}(M[a])}\right).$$

5.5. *Monotonicity properties of partial orders.* Suppose that $(a_1,\ldots,a_m, b_1,\ldots,b_n)$ is a uniformly distributed random permutation on $\mathfrak{S}_{m+n}$, the set of permutations on $m+n$ symbols. The permutation is to be viewed as the actual ranking of tennis skills of players $a_1,\ldots,a_m,b_1,\ldots,b_n$. We suppose that a player $x$ always loses to a player $y$ in a match if $x < y$. In a contest between teams $A = \{a_1,\ldots,a_m\}$ and $B = \{b_1,\ldots,b_n\}$, suppose that there is a partial order $\Theta$ between certain $a$'s and between certain $b$'s, for example, $a_1 < a_2$, $a_1 < a_3$, $b_2 < b_1$, ..., but no information about relative rankings between any $a_i$ and $b_j$. Such a situation arises if, to date, there have been numerous intrateam matches resulting, for example, in $a_1$ losing to $a_2$, $a_1$ losing to $a_3$, $b_2$ losing to $b_1$, and so on, but no interteam matches. We denote by $P(a_1 < b_1|\Theta)$ the conditional probability that $b_1$ defeats $a_1$ given the partial order $\Theta$.

Following numerous matches between members of the $A$ and $B$ teams, suppose the result has been a victory for $B$ in every case. Thus we now have information that the $a$'s have so far lost each match to the $b$'s. This induces a new partial order, $\Theta' = \Theta \cup \Theta''$, where $\Theta''$ consists of inequalities of the form $a_i < b_j$ for some collection of $i$ and $j$. We denote by $P(a_1 < b_1|\Theta')$ the conditional probability that $b_1$ defeats $a_1$ given the information in $\Theta'$.

Graham, Yao and Yao (1980) [see Graham (1982, 1983)] proved that

(5.10) $$P(a_1 < b_1|\Theta') \geq P(a_1 < b_1|\Theta),$$

and Shepp (1980) later gave another proof using the FKG inequality. As observed by Shepp, the additional knowledge with $\Theta'$ that a number of $a$'s severally have lost to a $B$-team member provides the basis for us to infer that the $A$-team is jointly inferior to the $B$-team. Therefore it is natural to expect a higher probability conditional on $\Theta'$, than conditional on $\Theta$, that $a_1$ loses to $b_1$ and (5.10) confirms this expectation.

Shepp (1980) constructed a suitable finite distributive lattice, a nonnegative MTP$_2$ measure $\mu$ and two decreasing indicator functions $f$ and $g$. Then he deduced (5.10) by an application of the FKG inequality. Along the same lines, we obtain generalizations of (5.10) using Shepp's MTP$_2$ measure by constructing three or more decreasing indicator functions $f_j$ and applying Theorem 1.1 or 1.3 to the functions $1 - f_j$. As an example, we state a result which follows from Theorem 1.1.

PROPOSITION 5.6. *Let $A_0$ be the subset of $\mathfrak{S}_{m+n}$ for which $A$ and $B$ have the complete order*:

$$A_0 = \{a_1 < \cdots < a_m\} \cap \{b_1 < \cdots < b_n\}.$$

*Suppose that $A_1, \ldots, A_4$ are subsets of $\mathfrak{S}_{m+n}$, each of which is an intersection of subsets of the form $a_i < b_j$. For any $A \subseteq \mathfrak{S}_{m+n}$, define*

$$\pi(A) := \frac{P(A_0 \cap A_4 \cap A)}{P(A_0 \cap A_4)}.$$

*Then*

$$\begin{aligned}
& 2\pi(A_1 \cap A_2 \cap A_3) \\
& \quad - [\pi(A_1 \cap A_2)\pi(A_3) + \pi(A_1 \cap A_3)\pi(A_2) + \pi(A_1)\pi(A_2 \cap A_3)] \\
& \quad + \pi(A_1)\pi(A_2)\pi(A_3) \\
& \quad - 2[\pi(A_1 \cap A_2) - \pi(A_1)\pi(A_2) + \pi(A_1 \cap A_3) \\
& \quad\quad - \pi(A_1)\pi(A_3) + \pi(A_2 \cap A_3) - \pi(A_2)\pi(A_3)] \leq 0.
\end{aligned}$$

By similar arguments, we can also generalize related results of Shepp (1980) and the XYZ conjecture [Shepp (1982)].

5.6. *Cumulants inequalities for probability distributions.* Our initial motivation for investigating the higher-order inequalities of FKG-type was to generalize numerous correlation inequalities well known for MTP$_2$ probability distributions in mathematical statistics; see Eaton (1987) and Karlin and Rinott (1980). We shall leave it to the reader to deduce from Theorems 1.1 and 1.3 bounds on the third-, fourth-, and fifth-order cumulants of those probability distributions.

The basis for much of Section 4 of Karlin and Rinott (1980) is a special case of the following result.

COROLLARY 5.7. *Let $\varphi$ and $\psi$ be MTP$_2$ functions on $\mathbb{R}^n$ and define*

$$\mathbb{E}(f) := \frac{\int \varphi(x)\psi(x)f(x)\,dx}{\int \varphi(x)\psi(x)\,dx}$$

*for any function $f$ for which the integrals converge. If $f_j$, $j = 1, 2, 3$, are nonnegative increasing functions on $\mathbb{R}^n$, then $\kappa_3'(f_1, f_2, f_3) \geq 0$.*

The proof of this result follows immediately from the generalized third-order FKG inequality once we note that the probability measure which is proportional to $\varphi(x)\psi(x)\,dx$, $x \in \mathbb{R}^n$, is an $\mathrm{MTP}_2$ measure. As applications of this result we can then deduce higher-order probability inequalities for any $\mathrm{MTP}_2$ random vector generalizing, for example, Example 4.1 of Karlin and Rinott (1980). In closing this section we provide, as a generalization of Proposition 4.1 of Karlin and Rinott (1980), a higher-order log-concavity property of exchangeable random variables.

PROPOSITION 5.8. *Let $X_1, \ldots, X_n$ be exchangeable random variables having a joint $\mathrm{MTP}_2$ probability density function $\varphi$. For $a \in \mathbb{R}$, define $c_0(a) = 1$ and*

$$c_m(a) = P(X_1 \leq a, \ldots, X_m \leq a),$$

$1 \leq m \leq n$. *Then*

(5.11)
$$2\frac{c_{m+2}(a)}{c_{m-1}(a)} - 3\frac{c_{m+1}(a)}{c_{m-1}(a)}\frac{c_m(a)}{c_{m-1}(a)}$$
$$+ \frac{c_m(a)^3}{c_{m-1}(a)^3} - 6\left[\frac{c_{m+1}(a)}{c_{m-1}(a)} - \frac{c_m(a)^2}{c_{m-1}(a)^2}\right] \leq 0,$$

$m = 1, \ldots, n-2$.

PROOF. For $x \in \mathbb{R}^n$, define

$$\psi(x) = \begin{cases} 1, & \text{if } x_1, \ldots, x_{m-1} \leq a, \\ 0, & \text{otherwise.} \end{cases}$$

It is well known that $\psi$ is $\mathrm{MTP}_2$. For $k = 1, 2, 3$, define

$$f_k(x) = \begin{cases} 1, & \text{if } x_{m+k-1} \leq a, \\ 0, & \text{otherwise.} \end{cases}$$

Then the functions $f_1$, $f_2$ and $f_3$ are decreasing. Moreover, $\mathbb{E}(f_k) = c_m(a)/c_{m-1}(a)$, $k = 1, 2, 3$; $\mathbb{E}(f_1 f_2) = \mathbb{E}(f_1 f_3) = \mathbb{E}(f_2 f_3) = c_{m+1}(a)/c_{m-1}(a)$; and $\mathbb{E}(f_1 f_2 f_3) = c_{m+2}(a)/c_{m-1}(a)$.

Now we apply Corollary 5.7 to the functions $1 - f_1$, $1 - f_2$ and $1 - f_3$. Simplifying the resulting expression, we obtain (5.11). □

By applying Proposition 5.8, we can obtain generalizations of other examples given by Karlin and Rinott [(1980), Section 4].

**6. Remarks on total positivity and inequalities of FKG-type.** In the development of inequalities of FKG-type, it will be instructive to study the case of indicator functions on $\mathbb{R}^2$. In what follows, for any $a \in \mathbb{R}$, we use the notation

$$I_a(t) = \begin{cases} 1, & t \geq a, \\ 0, & t < a, \end{cases}$$

for the indicator function of the interval $[a, \infty)$.

A prototypical increasing function $f$ on $\mathbb{R}^2$ is an indicator function of a "northeast" region, $[a, \infty) \times [b, \infty)$, so that $f$ is of the form

(6.1) $$f(u,v) = I_a(u)I_b(v) \equiv \begin{cases} 1, & u \geq a, v \geq b, \\ 0, & \text{otherwise,} \end{cases}$$

for some $a, b \in \mathbb{R}$. Let us consider the case of two such functions $f_j(u,v) = I_{a_j}(u)I_{b_j}(v)$, $(u,v) \in \mathbb{R}^2$, where $a_j, b_j \in \mathbb{R}$, $j = 1, 2$. In establishing the FKG inequality for these functions we may assume, by symmetry, that $a_1 \leq a_2$. We denote by $(X_1, X_2)$ the random vector corresponding to the density function $K$.

Suppose that $b_1 \leq b_2$. Then

$$Ef_1 = P(X_1 \geq a_1, X_2 \geq b_1) \leq 1$$

and

$$Ef_1 f_2 = Ef_2 = P(X_1 \geq a_2, X_2 \geq b_2) \geq 0.$$

Therefore

$$\text{Cov}(f_1, f_2) = Ef_1 f_2 - (Ef_1)(Ef_2) = Ef_2 - (Ef_1)(Ef_2)$$
$$= (Ef_2)(1 - Ef_1) \geq 0.$$

Next suppose that $b_1 > b_2$. Then

(6.2)
$$\text{Cov}(f_1, f_2) = Ef_1 f_2 - (Ef_1)(Ef_2)$$
$$= P(X_1 \geq a_2, X_2 \geq b_1)$$
$$\quad - P(X_1 \geq a_1, X_2 \geq b_1)P(X_1 \geq a_2, X_2 \geq b_2)$$
$$\equiv P(X_1 \geq a_1, X_2 \geq b_1)[1 - P(X_1 \geq a_1, X_2 \geq b_2)]$$
$$+ \begin{vmatrix} P(X_1 \geq a_1, X_2 \geq b_2) & P(X_1 \geq a_1, X_2 \geq b_1) \\ P(X_1 \geq a_2, X_2 \geq b_2) & P(X_1 \geq a_2, X_2 \geq b_1) \end{vmatrix}.$$

The first term in (6.2) clearly is nonnegative, so it remains to establish nonnegativity of the second term. To that end, we write

(6.3)
$$\begin{vmatrix} P(X_1 \geq a_1, X_2 \geq b_2) & P(X_1 \geq a_1, X_2 \geq b_1) \\ P(X_1 \geq a_2, X_2 \geq b_2) & P(X_1 \geq a_2, X_2 \geq b_1) \end{vmatrix}$$
$$= \begin{vmatrix} \int I_{a_1}(u)I_{b_2}(v)K(u,v)\,du\,dv & \int I_{a_1}(u)I_{b_1}(v)K(u,v)\,du\,dv \\ \int I_{a_2}(u)I_{b_2}(v)K(u,v)\,du\,dv & \int I_{a_2}(u)I_{b_1}(v)K(u,v)\,du\,dv \end{vmatrix}.$$

By two applications of the basic composition formula [Karlin (1968), page 17], once in $u$ and once in $v$, we see that (6.3) reduces to

$$
(6.4) \quad 4 \iint\limits_{u_1<u_2} \iint\limits_{v_1<v_2} \begin{vmatrix} I_{a_1}(u_1) & I_{a_1}(u_2) \\ I_{a_2}(u_1) & I_{a_2}(u_2) \end{vmatrix} \cdot \begin{vmatrix} I_{b_2}(v_1) & I_{b_2}(v_2) \\ I_{b_1}(v_1) & I_{b_1}(v_2) \end{vmatrix} \\ \times \begin{vmatrix} K(u_1,v_1) & K(u_1,v_2) \\ K(u_2,v_1) & K(u_2,v_2) \end{vmatrix} du_1\, dv_1\, du_2\, dv_2.
$$

The determinant $\det(I_{a_i}(x_j))$ is well known to be nonnegative if $x_1 < x_2$ and $a_1 < a_2$; indeed, the set of functions $\{I_{a_1}, I_{a_2} : a_1 < a_2\}$ is an example of a *weak Tchebycheff system* [cf. Karlin (1968), Chapter 1]. Therefore the first two determinants in (6.4) are nonnegative on the region $\{u_1 < u_2, v_1 < v_2\}$. Since $K$ is $TP_2$, then the third determinant also is nonnegative on the same region. Hence $\text{Cov}(f_1(X_1, X_2), f_2(X_1, X_2)) \geq 0$.

Next suppose that $f_1$ and $f_2$ are functions on $\mathbb{R}^n$ of the form

$$(6.5) \qquad f_j(x_1,\ldots,x_n) = I_{a_{j,1}}(x_1) I_{a_{j,2}}(x_2) \cdots I_{a_{j,n}}(x_n),$$

for $j = 1, 2$. We assume, by induction, that the FKG inequality holds for all functions of the above form for all dimensions up to $n-1$. Thus, we have

$$(6.6) \qquad f_j(x_1,\ldots,x_n) = g_j(x_1,\ldots,x_{n-1}) I_{a_{j,n}}(x_n),$$

where $g_j(x_1,\ldots,x_{n-1}) = I_{a_{j,1}}(x_1) \cdots I_{a_{j,n-1}}(x_{n-1})$, $j = 1, 2$, and

$$\text{Cov}(g_1(X_1,\ldots,X_{n-1}), g_2(X_1,\ldots,X_{n-1})) \geq 0.$$

Using the standard method based upon conditional expectations, we now complete the inductive step; we underscore that this inductive step is well known and we provide details only for the sake of completeness. Denoting by $E_{X_1,\ldots,X_{n-1}|X_n}$ expectation with respect to the probability distribution of the random variables $X_1,\ldots,X_{n-1}$ conditional on $X_n$, it follows from (6.6) and the law of total probability that

$$
(6.7) \begin{aligned} Ef_1 f_2 &= E_{X_n} E_{X_1,\ldots,X_{n-1}|X_n} f_1 f_2 \\ &= E_{X_n} I_{a_{1,n}}(X_n) I_{a_{2,n}}(X_n) E_{X_1,\ldots,X_{n-1}|X_n} \prod_{j=1}^{2} g_j(X_1,\ldots,X_{n-1}). \end{aligned}
$$

Since $(X_1,\ldots,X_n)$ is $MTP_2$, then $(X_1,\ldots,X_{n-1})|X_n$ is also $MTP_2$ [see Sarkar (1969) and Karlin and Rinott (1980)]. Therefore, by inductive hypothesis,

$$E_{X_1,\ldots,X_{n-1}|X_n} g_1 g_2 \geq (E_{X_1,\ldots,X_{n-1}|X_n} g_1)(E_{X_1,\ldots,X_{n-1}|X_n} g_2),$$

and then it follows from (6.7) that

$$Ef_1 f_2 \geq E_{X_n} \psi_1(X_n) \psi_2(X_n),$$

where

$$\psi_j(x_n) = I_{a_{j,n}}(x_n) E_{X_1,\ldots,X_{n-1}|X_n=x_n} g_j(X_1,\ldots,X_{n-1}),$$

$j = 1, 2$. Since $(X_1, \ldots, X_n)|X_n$ is $\text{MTP}_2$ and $g_j$ is increasing, then [cf. Sarkar (1969) and Karlin and Rinott (1980), page 484, Theorem 4.1] the function $E_{X_1,\ldots,X_{n-1}|X_n=x_n} g_j(X_1,\ldots,X_{n-1})$ is increasing in $x_n$; hence $\psi_1$ and $\psi_2$ are both increasing, so that $E\psi_1\psi_2 \geq E\psi_1 E\psi_2$. Therefore we obtain

$$Ef_1 f_2 \geq E\psi_1\psi_2 \geq E\psi_1 E\psi_2 \equiv Ef_1 Ef_2,$$

and the proof of the FKG inequality for functions of the type (6.5) is complete.

Having established the FKG inequality (1.2) for all indicator functions of the form (6.5), we observe that the functional $(f_1, f_2) \to \text{Cov}(f_1, f_2)$ is bilinear. Hence (1.2) holds for all functions $f_1$ and $f_2$ on $\mathbb{R}^n$ of the form

$$f(x_1,\ldots,x_n) = \sum_{i=1}^{r} c_i \prod_{j=1}^{n} I_{a_{i,j}}(x_j),$$

$(x_1,\ldots,x_n) \in \mathbb{R}^n$, where $r \in \mathbb{N}$ and $c_i \geq 0$, $i = 1,\ldots,r$. Denoting by $\mathcal{M}_+(\mathbb{R}^n)$ the set of positive Borel measures on $\mathbb{R}^2$, it then follows by approximation arguments that the FKG inequality (1.2) holds for all functions $f_1$ and $f_2$ of the form

(6.8) $$f(x_1,\ldots,x_n) = \int_{\mathbb{R}^n} \prod_{j=1}^{n} I_{a_j}(x_j) \, d\nu(a_1,\ldots,a_n),$$

$(x_1,\ldots,x_n) \in \mathbb{R}^n$, where $\nu \in \mathcal{M}_+(\mathbb{R}^n)$.

It is well known that the set of functions of type (6.8) is a proper subset of the class of increasing functions on $\mathbb{R}^n$ and contains all cumulative distribution functions on $\mathbb{R}^n$.

Turning to the third-order FKG inequality in Theorem 1.1, we can also establish that result for the class of functions (6.1). Suppose that we have three indicator functions, $f_j(u,v) = I_{a_j}(u) I_{b_j}(v)$, $(u,v) \in \mathbb{R}^2$, $j = 1,2,3$. To establish Theorem 1.1 for $f_1$, $f_2$ and $f_3$, we may assume, by symmetry, that $a_1 \leq a_2 \leq a_3$. Then the proof requires that we resolve six cases, each corresponding to an ordering of $b_1$, $b_2$ and $b_3$. In what follows, we shall denote $P(X_1 \geq a_i, X_2 \geq b_j)$ by $\rho_{ij}$.

CASE 1 $b_1 \leq b_2 \leq b_3$. In this case, $f_i f_j \equiv f_j$ for $i \leq j$. Therefore

$$\kappa'_3 = 2\rho_{33} - [\rho_{22}\rho_{33} + \rho_{33}\rho_{22} + \rho_{11}\rho_{33}] + \rho_{11}\rho_{22}\rho_{33}$$
$$= (2 - \rho_{11})(1 - \rho_{22})\rho_{33},$$

which, clearly, is nonnegative.

CASE 2 $b_1 \leq b_3 \leq b_2$. Here, we have
$$\kappa'_3 = 2\rho_{32} - [\rho_{22}\rho_{33} + \rho_{33}\rho_{22} + \rho_{11}\rho_{32}] + \rho_{11}\rho_{22}\rho_{33}$$
$$= (2 - \rho_{11})(\rho_{32} - \rho_{33}\rho_{22}).$$

By the FKG inequality, $\rho_{32} - \rho_{33}\rho_{22} = \mathbb{E}(f_3 f_2) - \mathbb{E}(f_3)\mathbb{E}(f_2) \geq 0$; therefore $\kappa'_3 \geq 0$.

CASE 3 $b_2 \leq b_1 \leq b_3$. In this case, we have
$$\kappa'_3 = 2\rho_{33} - [\rho_{21}\rho_{33} + \rho_{33}\rho_{22} + \rho_{11}\rho_{33}] + \rho_{11}\rho_{22}\rho_{33}$$
$$= \rho_{33}[1 - \rho_{21} + (1 - \rho_{11})(1 - \rho_{22})],$$

which, clearly, is nonnegative.

CASE 4 $b_2 \leq b_3 \leq b_1$. Here,
$$\kappa'_3 = 2\rho_{31} - [\rho_{21}\rho_{33} + \rho_{31}\rho_{22} + \rho_{11}\rho_{33}] + \rho_{11}\rho_{22}\rho_{33}$$
$$= (1 - \rho_{22})(\rho_{31} - \rho_{33}\rho_{11}) + \rho_{31} - \rho_{21}\rho_{33}.$$

By the FKG inequality, both $\rho_{31} - \rho_{33}\rho_{11} = \mathbb{E}(f_3 f_1) - \mathbb{E}(f_3)\mathbb{E}(f_1)$ and $\rho_{31} - \rho_{33}\rho_{21} = \mathbb{E}(f_3 f_2 f_1) - \mathbb{E}(f_3)\mathbb{E}(f_2 f_1)$ are nonnegative. Therefore $\kappa'_3 \geq 0$.

CASE 5 $b_3 \leq b_1 \leq b_2$. In this case, we have
$$\kappa'_3 = 2\rho_{32} - [\rho_{22}\rho_{33} + \rho_{31}\rho_{22} + \rho_{11}\rho_{32}] + \rho_{11}\rho_{22}\rho_{33}$$
$$= (1 - \rho_{11})(\rho_{32} - \rho_{33}\rho_{22}) + \rho_{32} - \rho_{31}\rho_{22}.$$

By the FKG inequality, both $\rho_{32} - \rho_{33}\rho_{22} = \mathbb{E}(f_3 f_2) - \mathbb{E}(f_3)\mathbb{E}(f_2)$ and $\rho_{32} - \rho_{31}\rho_{22} = \mathbb{E}(f_1 f_2 f_3) - \mathbb{E}(f_3 f_1)\mathbb{E}(f_2)$ are nonnegative. Therefore $\kappa'_3 \geq 0$.

CASE 6 $b_3 \leq b_2 \leq b_1$. In this case, we have
$$\kappa'_3 = 2\rho_{31} - [\rho_{21}\rho_{33} + \rho_{31}\rho_{22} + \rho_{11}\rho_{32}] + \rho_{11}\rho_{22}\rho_{33}$$
$$= (1 - \rho_{22})(\rho_{31} - \rho_{33}\rho_{11}) + \rho_{33}(1 - \rho_{21}) + \rho_{11}(\rho_{33} - \rho_{32}).$$

By the FKG inequality, $\rho_{31} - \rho_{33}\rho_{11} = \mathbb{E}(f_3 f_1) - \mathbb{E}(f_3)\mathbb{E}(f_1) \geq 0$. Also, since $1 - f_2 \geq 0$, $\rho_{33} - \rho_{32} = \mathbb{E}(f_3) - \mathbb{E}(f_2 f_3) = \mathbb{E} f_3(1 - f_2) \geq 0$. Therefore $\kappa'_3 \geq 0$.

Having resolved the case in which the $f_j$ are indicator functions, we apply the multilinearity of the functional $(f_1, f_2, f_3) \to \kappa'_3(f_1, f_2, f_3)$, and an approximation argument, to deduce nonnegativity of $\kappa'_3(f_1, f_2, f_3)$ for all $f_j$ of the class (6.8) with $n = 2$.

In the case of $\kappa'_4$, we have carried out the case-by-case analysis (consisting of 24 cases); as regards $\kappa'_5$, we unhesitatingly entrust the analysis (of all 120 cases) to the reader.

Even in the case of the classical FKG inequality, the method of indicator functions appears to be new. The technique has the obvious drawback that it yields the classical FKG inequality only for a proper subset of the class of all increasing functions, and that too by way of a case-by-case analysis. However, the method has the advantage that it points the way toward generalizations of that inequality; indeed, the method of indicator functions is the means by which we first discovered instances of the nonnegativity of $\kappa_3'$.

In developing the method of indicator functions, the appearance of the Binet–Cauchy formula is noteworthy. To explain, we first remark that in a previous paper, Gross and Richards (1995) developed a Binet–Cauchy formula in the context of finite reflection groups and obtained an analog of the FKG inequality for one particular reflection group. We speculate that analogs of those inequalities exist in any context in which Binet–Cauchy formulas are available, including those of Karlin and Rinott (1988).

A general context in which formulas of Binet–Cauchy type arise is given by Selberg (1956) in a fundamental paper laying the groundwork for the general theory of what is now known as *Selberg's trace formula*. Let $S$ denote a Riemannian space, with local coordinates $t^1, \ldots, t^n$ and positive-definite metric

$$ds^2 = \sum g_{ij}\, dt^i\, dt^j,$$

where the functions $g_{ij}$ are analytic in $t^1, \ldots, t^n$. We assume that we have a locally compact group $G$ of isometries of $S$, with a transitive action. A function $K: S \times S \to \mathbb{C}$ is *point-pair invariant* if $K(gx, gy) = K(x, y)$ for all $x, y \in S$ and $g \in G$. One of the basic problems considered by Selberg (1956) is the analysis of the spectral theory of the integral operator

$$f \mapsto \int_S K(x,y) f(y)\, dy,$$

where $dy$ denotes the invariant measure on $S$ derived from the metric $ds^2$ and $f$ is a suitable function.

Let $\Gamma$ denote a discrete subgroup of $G$ which acts properly and discontinuously on $S$. Denote by $\mathfrak{D}$ the fundamental domain of $\Gamma$. With $\mathrm{U}(d)$ denoting the group of unitary $d \times d$ matrices, let $\chi: \Gamma \to \mathrm{U}(d)$ be a unitary representation of $\Gamma$. Let $F$ be a complex, vector-valued function on $S$ such that, for all $x \in S$ and $M \in \Gamma$, $F(Mx) = \chi(M) F(x)$. Then it is a readily established, but important, result that

$$\int_S K(x,y) F(y)\, dy = \int_\mathfrak{D} K_\chi(x,y) F(y)\, dy$$

where

(6.9) $$K_\chi(x,y) := \sum_{M \in \Gamma} \chi(M) K(x, My);$$

36 D. ST. P. RICHARDSsee Selberg [(1956), page 59]. Formally, there exists an eigenfunction expansion of $K_\chi$:

$$K_\chi(x,y) = \sum c_i F_i(x) F_i(y)^*,$$

where the $c_i$ are constants, the $F_i$ are eigenfunctions of a class of differential operators invariant under the action of $G$, and $F_i(y)^*$ is the transpose of the complex conjugate of $F_i(y)$. What is of interest to us here is that the function $K_\chi$ possesses properties of the determinant function [defined by Gross and Richards (1995)] for finite reflection groups. Indeed, it follows from (6.9) that

$$K_\chi(x, My) = K_\chi(x, y) \chi(M^{-1})$$

for all $x, y \in S$, $M \in \Gamma$; this property generalizes the familiar sign-change behavior of the classical determinants under the interchange of rows or columns. The function $K_\chi$ also satisfies a generalized Binet–Cauchy formula: If $K$ and $L$ are point-pair invariant functions, and

$$Q(x,y) := \int_S K(x,t) L(t,y)\, dt$$

then $Q$ is point-pair invariant and

$$Q_\chi(x,y) = \int_\mathfrak{D} K_\chi(x,t) L_\chi(t,y)\, dt.$$

We conjecture that inequalities of FKG-type exist in this context. Owing to the potential for applications in multivariate statistical analysis, of special interest for us will be the case in which $S$ is the space of positive-definite symmetric matrices.

## REFERENCES

Ahlswede, R. and Daykin, D. E. (1978). An inequality for the weights of two families of sets, their unions and intersections. *Z. Wahrsch. Verw. Gebiete* **43** 183–185. MR491189

Barrett, W. W. and Feinsilver, P. J. (1981). Inverses of banded matrices. *Linear Algebra Appl.* **41** 111–130. MR649720

Batty, C. J. K. and Bollmann, H. W. (1980). Generalised Holley–Preston inequalities on measure spaces and their products. *Z. Wahrsch. Verw. Gebiete* **53** 157–173. MR580910

Cameron, P. J. (1987). On the structure of a random sum-free set. *Probab. Theory Related Fields* **76** 523–531. MR917677

Cartier, P. (1974). Inégalités de corrélation en mécanique statistique. *Lecture Notes in Math.* **383** 1–33. Springer, Berlin. MR418714

Chen, M.-F. and Wang, F.-Y. (1993). On order-preservation and positive correlations for multidimensional diffusion processes. *Probab. Theory Related Fields* **95** 421–428. MR1213199

Eaton, M. L. (1987). *Lectures on Topics in Probability Inequalities. CWI Tract* **35**. Stichting Mathematisch Centrum, Centrum voor Wiskunde en Informatica, Amsterdam. MR889671

Department of Statistics
Pennsylvania State University
University Park, Pennsylvania 16802
USA
e-mail: richards@stat.psu.edu